\documentclass[preprint,12pt]{elsarticle}

\usepackage{amssymb}
\usepackage{amsmath}
\usepackage[a4paper,margin=1in]{geometry}
\usepackage{graphicx}
\usepackage[colorlinks=true, allcolors=blue]{hyperref}
\usepackage{longtable}
\usepackage{makecell}
\usepackage{setspace}
\usepackage{subcaption}
\usepackage{times}
\usepackage{xurl}

\journal{Applied Energy}

\begin{document}

\begin{frontmatter}



\title{\textbf{Integrating Data Center Flexibility into Power System Operations: A Framework for Temporal Redispatch}}

\author[orie]{Huansong Li \corref{cor1}}
\ead{huansong@utexas.edu}

\author[orie]{Erhan Kutanoglu}
\ead{erhank@austin.utexas.edu}

\author[che]{Michael Baldea}
\ead{mbaldea@che.utexas.edu}

\author[che]{Ilias Mitrai}
\ead{imitrai@che.utexas.edu}

\cortext[cor1]{Corresponding author. Tel.: +1 404-547-2059.}

\affiliation[orie]{organization={Operations Research and Industrial Engineering, The University of Texas at Austin},
                   addressline={204 E. Dean Keeton Street, C2200}, 
                   city={Austin},
                   postcode={78712-1591}, 
                   state={Texas},
                   country={United States}}

\affiliation[che]{organization={McKetta Department of Chemical Engineering, The University of Texas at Austin},
                  addressline={200 E. Dean Keeton St. Stop C0400}, 
                  city={Austin},
                  postcode={78712-1589}, 
                  state={Texas},
                  country={United States}}

\begin{abstract}
Rapid growth in data center electricity demand, driven by artificial intelligence (AI) training and inference, is expected to increase stress conditions for power grids. Operational flexibility from data centers offers a near-term lever to relieve such stresses. Existing models of data center flexibility tend to bundle economic incentives with physical constraints, obscuring the tradeoff between economic willingness to shift load and technical limitations to do so. This paper proposes a temporal-redispatch focused framework that decouples data center demand flexibility into two independently parameterized dimensions: (1) systematic flexibility, capturing the economic environment through penalties for temporal shifting of data centers' scheduled demand, and (2) operational flexibility, capturing the physical envelope through the demand shifting time window and the limit on shifted data center demand. The framework relies on a direct current optimal power flow (DC-OPF) linear program which enables annual parametric analysis at hourly resolution. The framework is studied using a synthetic grid mimicking the Texas transmission system populated with the 2026 Adjusted Forecast from ERCOT (Electric Reliability Council of Texas). Results show that operational flexibility imposes a physical lower bound on unsatisfied demand that no economic incentive can overcome; longer shifting windows and higher limits on shifted data center demand reduce operational cost with different impacts on grid utilization; and, temporal redispatch shaves system-wide net load peaks while partially compensating for localized transmission bottlenecks that emerge under modeled coal-retirement scenarios.
\end{abstract}


\begin{highlights}
\item Framework separates data-center economic willingness from technical capability
\item Annual hourly analyses in DC-OPF linear program on a synthetic ERCOT-like grid
\item Data-center shifting constraints impose a lower bound on unmet demand
\item Longer windows and higher scheduled-demand caps expand feasible redispatch
\item Local coal retirements can propagate reliability impacts through congestion
\end{highlights}

\begin{keyword}

Data centers \sep
Temporal redispatch \sep
DC-OPF \sep
Power system flexibility \sep
Transmission congestion \sep
Demand response

\end{keyword}

\end{frontmatter}


\section{Introduction}
\label{sec: introduction}

The global energy landscape is being reshaped by an unprecedented increase in computational demand. Artificial intelligence (AI) training and inference, large-scale cloud services, and continued growth in data processing have driven a step change in data center deployment and energy consumption. According to the International Energy Agency (IEA), global data center electricity consumption reached 415 terawatt-hours (TWh) in 2024 — a figure that already exceeds the total annual electricity consumption of many developed nations — and is projected to more than double to 945 TWh by 2030 \cite{IEA20251}.

The impact of data center proliferation on power grids within the United States is intense. A 2024 report by Lawrence Berkeley National Laboratory (LBNL) established the baseline statistics: data centers consumed 176 TWh of electricity in 2023, approximately 4.4\% of total U.S. electricity consumption, with base-case projections placing this share at 6.7\% to 12.0\% (325 TWh to 580 TWh) by 2028 \cite{LBNLP1WC7Q}. The growth is geographically concentrated. Northern Virginia currently hosts the largest data center cluster in the world \cite{NAP29101}, and grid capacity exhaustion in established hubs has begun to redirect new development toward deregulated energy markets capable of supporting large interconnections \cite{IEA20251}.

The Electric Reliability Council of Texas (ERCOT), which manages the power grid for 90\% of the state's electrical load, has recently emerged as a center for new data center load growth \cite{EIA2024}. ERCOT's unique market structure, isolated interconnection, and abundant wind and solar resources attract new large load projects. However, the sheer volume of incoming interconnection requests is a challenge for the grid. The ERCOT 2025 Long-Term Load Forecast (LTLF) projects that the summer peak demand will reach 218 gigawatts (GW) by 2031 \cite{ERCOT20251}. This represents a massive upward revision from the 2024 LTLF, which had projected a peak demand of only 150 GW for the same time frame \cite{ERCOT20252}. This increase is almost entirely attributable to projected data center load growth, formally attested to by Transmission Service Provider (TSP) Officer Letters \cite{ERCOT20252}. Specifically, the 2025 TSP-provided forecast indicates that data center growth for the year 2030 is revised from a projected 29.6 GW (in the 2024 forecast) to 77.9 GW (in the 2025 forecast) \cite{ERCOT20251}. To put this in perspective, the projected five-year load growth currently sitting in ERCOT interconnection queues is equivalent to adding multiple times the peak electrical load of major metropolitan areas to the Texas grid in less than a decade \cite{GS2025, ARXIV260529053}.

Because a single hyperscale facility can draw the output of a small-to-medium conventional power plant, this magnitude and pace of integration pose direct risks to grid stability, including localized congestion, voltage depression, and line thermal overload \cite{EIA2024}. The structural problem is the mismatch between the pace at which data center load commitments materialize and the pace at which bulk-system infrastructure can be expanded. New generation and transmission buildout requires planning, permitting, and construction horizon spanning multiple years to a decade. This horizon is getting longer: an IEA survey of industry participants reports that procurement lead times for large power transformers have almost doubled since 2021, now reaching up to four years, with high-voltage direct-current cables exceeding five years \cite{IEA20252}. Hyperscale data center commitments enter interconnection queues and their loads are reflected in system forecasts well before the corresponding bulk-system expansions can be delivered. Data center developers can phase deployment or co-locate behind-the-meter generation to provide power ahead of network upgrades. Increased electricity demand from early data center investments is being realized on the grid long before the underlying generation and transmission investments can be made \cite{ACM20241}. Because these constraints affect the whole system, they slow the grid and generation expansion, reinforcing the case for operational flexibility of data center demand shifting as a near-term solution.

There is one more factor to consider for data center demand flexibility. Recent studies \cite{PE2024, EP4GKW61} have noted that as the generation mix shifts toward intermittent resources such as wind and solar, the economic cost of transmission congestion has risen significantly, necessitating more responsive demand-side participation. Consequently, the reliability of the existing grid must be maintained through advanced operational strategies rather than relying solely on physical capacity expansion. On one hand, a data center requires continuous electricity to maintain server uptime, prevent data corruption, and satisfy stringent Service Level Agreements (SLAs) \cite{NI2025}. On the other hand, the digital workload a data center processes, such as batch computing and AI model training, can be shifted across spatial or temporal clusters. Data centers can rapidly modulate their power loads whenever the grid's other needs take priority, making them a responsive candidate for advanced demand response and grid ancillary services \cite{SEIR2025}. The inherent agility of data centers positions them as a critical source of flexibility for power systems that manage high penetrations of intermittent renewable energy. However, properly harnessing this flexibility requires a precise definition and parameterization of that  ``flexibility.'' Existing literature does not distinguish between what data centers are \textit{technically capable} of doing and what the grid is \textit{economically willing} to consider for demand shifting under power flow constraints. Defining and parameterizing flexibility is a critical gap in grid optimization research, providing the primary motivation for this paper.

In this paper, we introduce a novel framework for data center load redispatch that strictly categorizes flexibility into two hierarchical mathematical dimensions: 
\begin{itemize}
    \item Systematic Flexibility: defined as the system-wide economic ``willingness'' of the grid operator to accept costs associated with unserved energy or schedule deviations. In the context of mathematical optimization, systematic flexibility is dictated entirely by the penalty structures embedded within the objective function of the power system's operation. Specifically, it represents the macro-level economic environment, governed by the penalty parameter applied to unsatisfied grid demand, and the penalty parameter applied to data center demand shifting. These economic weights determine how aggressively the optimization model will attempt to balance standard generation costs with the financial repercussions of interrupting or deferring data center loads.
    \item Operational Flexibility: defined as the physical constraints within which the model can actually maneuver. It defines the ``feasible region'' of physical actions that the data center can technically execute without violating its core computational demands or thermal limits. In the proposed mathematical framework, operational flexibility is enforced strictly through physical constraints, defined by parameters such as the time window — which dictates the maximum duration a load can be deferred or advanced before it must be serviced — and the limit parameter that caps the magnitude of allowable load addition relative to the baseline data center demand.
\end{itemize}

This separation between the two types of flexibility and their parameterization distinguishes economic incentives from the technical capabilities typically seen in the literature. By formally modeling systematic and operational flexibility dimensions into a unified temporal redispatch model, this paper provides a mathematical optimization-based approach. Rather than modeling the data center as a generic, simplistic load-shedding resource, the proposed framework allows grid operators to evaluate the relationship and tradeoffs between economics driven penalties and physical or technological load-shifting capabilities. Within the scope of this paper, redispatch and shift are two interchangeable terms which refer to the moving (or shifting) of data center load from its scheduled time period (hour) to another time period within a limited (and parameterized) time window.

The remainder of this paper is structured as follows: Section \ref{sec: literature} provides a comprehensive literature review, focusing on the definition of flexibility from different dimensions and recent papers in the optimization field that investigate power system flexibility, joint capacity expansion, and spatiotemporal load shifting. Section \ref{sec: model} presents the mathematical formulation of the temporal redispatch model, detailing the objective function, physical network constraints, and the integration of the systematic and operational flexibility parameters. Section \ref{sec: experiment} provides a detailed experiment design with a series of case studies that apply the proposed model to a constrained power grid topology with coal generators reserved or partially retired, analyzing the experimental settings and the resulting impacts on load fulfillment and generation variability. Finally, Section \ref{sec: conclusion} offers a conclusion, summarizing the implications of data center flexibility and further opportunities for future research.

\section{Literature Review}
\label{sec: literature}

The power systems literature, national laboratories, and regulatory bodies have long recognized that demand-side integration involves two distinct dimensions: the economic incentives that alter consumer's behavior, and the physical or technical capabilities of the underlying assets. Foundational reports from the early 2010s explicitly established this dichotomy \cite{PEC2010, EPRI1016987, LLNL1880939, IEA2018, IEA2019, IEA2014, SU17177891, NYSERDA2014, APP9183796}. In these frameworks,``technical flexibility'' was defined as the physical ability of a resource to repeatedly change its operational behavior, governed by mechanical ramp rates, minimum up/down times, and storage capacities. ``Economic flexibility'' was defined as the willingness of a market participant to actually execute those changes in response to appropriate financial incentives, dynamic pricing signals, or penalty avoidance \cite{PEC2010}. For instance, a power plant might possess the technical capability to ramp down quickly, but lack the economic flexibility to do so if bound by rigid ``take-or-pay'' fuel contracts \cite{IEA2018}. As the power sector transitioned to smart grids and demand response, this paradigm was applied to large industrial consumers. A recent study indicated that although facilities possessed significant technical capabilities for load shedding, they historically exhibited low participation rates due to misaligned operational priorities and insufficient economic incentives \cite{NI2025}. This historical understanding serves as the theoretical precursor to the systematic and operational flexibility dimensions introduced in our framework.

The traditional approach to power system planning typically treats electricity demand as an exogenous, inflexible input and rely on mixed-integer linear programming (MILP) to optimize the addition or expansion of generators and transmission lines \cite{IISET2026, OE2019, EJOR202106024, MP2013, RSTA20190436}. For example, a multi-stage stochastic MILP for electric power infrastructure that accounts for the renewable intermittency and unit-commitment constraints was formulated in \cite{OE2019}. A later study proposed a MILP framework for jointly optimizing generation and transmission expansion that handles realistic discrete investment options and operating reserves \cite{EJOR202106024}. These frameworks remain effective for grid assets but do not account for the magnitude, deployment speed, or load-shifting capabilities of data centers. Recognizing this critical gap, more recent optimization work has shifted toward the \textit{joint planning} of power systems and data centers. A recent study proposed a joint planning framework applied to PJM interconnection that co-optimizes the long-term development of data centers and power systems over a 15-year horizon \cite{IJEPES2025111377}. In a complementary direction, another study used a near-optimal generation-portfolio approach on a Texas Interconnection model to explore how cost-effective expansion responds to large-load growth from data centers, cryptocurrency mining, and industrial loads \cite{NPJCE2025}. Another study examined the cost–emission tradeoffs of flexible data centers, showing that increased flexibility can lower costs but may shift emissions depending on the generation mix \cite{W34065}. The investment-economics literature continues to evolve: Recently, a study proposed a stochastic-programming framework for optimal energy portfolios of data centers under market uncertainty \cite{AE2026}. Although highly effective for strategic decisions, these frameworks operate at investment timescales and do not resolve the hour-to-hour interaction between an existing transmission topology, a constrained generation mix, power flow physics, and dynamically dispatchable data center loads, which is the combination of factors we seek to explore.

Because digital workloads can be routed across geographically separated facilities and deferred within a given facility, much of the data center flexibility literature couples spatial migration with temporal deferral. A substantial strand of this \textit{spatiotemporal} work operates at the distribution-network or data center-microgrid level, where voltage regulation, feeder thermal limits, and reactive-power management are key considerations and the network physics calls for specific formulations rather than transmission-level power flow. A robust mixed-integer nonlinear programming (MINLP) model was developed for the spatiotemporal load management of data centers on a modified IEEE 33-bus active distribution network, framing the data center as a ``virtual battery'' that mitigates renewable curtailment through task migration across data centers and workload deferral within each data center \cite{RSER2024114793}. An emission-aware resource-planning scheme for data center microgrids that jointly considers batch workload scheduling and risk management was formulated in \cite{IEEETIA2018}. A distributed real-time energy-management algorithm for data center microgrids based on consensus-style coordination was proposed in \cite{IEEETSG2018}. A game-theoretic framework for data center demand response in deregulated electricity markets was developed in \cite{IEEETSG2019}, focusing on the distribution-level interaction between data centers and a local utility. A comprehensive survey in \cite{RSER2019109466} catalogued load-modeling, load-regulation, and incentive-design strategies developed for data center participation in demand response. A more recent review in \cite{AE2021} surveyed the state-of-the-art for integrated data center and smart-grid systems. These distribution-level studies provide models of physical and economic interaction, but their findings are scoped to a single feeder or a microgrid and typically rely on the assumption that geographically distributed data centers are operated by a single entity capable of physically migrating computational workloads between sites, which do not transfer naturally to transmission systems or to single-zone independent system operator (ISO) operations.

A parallel strand of spatiotemporal work operates at the transmission level, capturing the interaction among bulk-system congestion, generation dispatch, and data center load-shifting flexibility. The ``virtual link'' abstraction, a non-physical pathway that captures the transfer of computing-driven power loads across space and time, was introduced by \cite{EPSR2020106723} within a market-clearing framework. This abstraction was subsequently extended along two distinct directions: one study established the economic properties of a market-clearing formulation that remunerates data centers for space–time load-shifting flexibility based on nodal price differences across space and time \cite{AE2022}. Another study later developed a computational framework that quantifies the operational flexibility associated with such shifting in electricity markets \cite{CCE2023}. A recent study quantified the spatial-flexibility benefit on a modified IEEE 73-bus system by co-optimizing generation dispatch, security reserves, and flexible data center loads \cite{ARXIV251108759}. The work found that enabling spatial redistribution to solar-rich regions eliminated all transmission violations, lowering peak line utilization and solar curtailment. Another study embedded data centers into a distributed economic dispatch problem via a G/G/1 queueing-theoretic representation paired with a Quality-of-Service(QoS)–based cost function \cite{TSUSC20253538817}. The work decomposed the centralized problem into a tractable Lagrangian dual and solved it on a 5-bus PJM test system, demonstrating that optimized server utilization can alleviate network congestion while respecting transmission constraints. Two methodological observations emerge from this research strand. First, although the test systems are nominally at the transmission level, they are predominantly abstract IEEE benchmarks (the 5-bus PJM, IEEE 30-bus, IEEE 73-bus test grids) rather than regional interconnect-scale topologies calibrated to a specific ISO's projected load profile. Recent work in \cite{ARXIV251023586} on spatial-aggregation impacts in capacity expansion planning demonstrated that the choice between zonal and nodal resolution can distort investment and operational decisions, motivating analyses on realistic large-scale networks. Second, most existing formulations embed economic penalties directly within physical or market-clearing constraints, intermixing what a data center is technically capable of doing with what the market is willing to pay. This entanglement is convenient for compact formulations but obscures the tradeoff between economic incentives  and physical capabilities of data centers, which we seek to address in this paper.

The literature reviewed above provides rich modeling machinery for data center–grid interaction, but the spatial dimension may carry assumptions that are not universally applicable. Multi-site coordination requires that these sites share an operator capable of redispatching workloads \cite{IEEETSC2020, IEEETSG2016} but cross-border data-sovereignty rules constrain such routing in many jurisdictions and end-user latency requirements limit workload migration distances. In single-zone congestion contexts — common inside large ISOs, where data center capacity within a constrained transmission zone is local and no compatible peer site at the computational scale exists — temporal redispatch becomes the main tool for providing grid flexibility. Time-based flexibility is therefore best understood as its own problem class rather than a special case of spatiotemporal flexibility, and it has been studied across a hierarchy of timescales. At the operational (i.e., hourly) timescale at which our framework sits, the literature is comparatively sparse: one study examined data center participation in demand response through an adaptive policy that enforces workload QoS guarantees while shedding power-grid-coincident load — a single-facility temporal-only treatment that demonstrates the operational feasibility of bulk-load deferral \cite{TSUSC20213077254}. Nevertheless, the work does not model the underlying transmission network or the projected large-load growth, which motivates our paper. In sub-hourly resolution, a study modeled data centers providing frequency-regulation services by approximating server power consumption as a piece-wise bi-linear function of CPU frequency and utilization rate, embedded within an optimization framework that determines bidding strategies in coupled energy and regulation markets \cite{AE2019}. Similarly, a later study proposed a data-driven linearization framework to co-optimize generation schedules and fast frequency response (FFR) from data centers \cite{ARXIV251001050}. Their quantitative assessment introduced a  ``Marginal Flexibility Value'' metric to quantify the diminishing returns from increasing fast-responding data center share. Previous work in \cite{AE2018}  characterized the dynamic flexibility provided by data center uninterruptible power supply systems and battery energy storage, demonstrating their viability for frequency containment without interrupting computing workloads. A complementary recent direction is carbon-aware computing: one study developed probabilistic performance guarantees for data center scheduling that tracks grid emissions signals over time, providing a temporal-only formulation focusing on emissions rather than congestion \cite{IEEETPS2025}. These studies define a timescale hierarchy, from sub-second FFR through sub-hourly frequency regulation to hourly redispatch, latter of  which the framework developed in this paper sits, most relevant to ISO operational planning and most near-term solution to address the rapid data center interconnection growth.

In short, three recurring features cut against direct operational usability of data center flexibility under projected AI-driven large-load growth. First, the majority of the work models spatial and temporal load-shifting jointly, despite the operational reality that spatial migration may not be available in all conditions; temporal-only redispatch at the hourly operational timescale is therefore understudied relative to its operational importance for ISOs managing localized transmission stress. Second, previous formulations tend to bundle economic incentives within physical constraints, so that questions about how aggressively to incentivize shifting cannot be cleanly separated from questions about what data centers are physically capable of doing. Third, transmission-level analyses of data center flexibility have largely been carried out on abstract IEEE test grids or compressed temporal horizons rather than on regional interconnect-scale topologies under realistic projected load. Even though temporal redispatch acts on a localized resource (say, a bus), the congestion it relieves is determined by system-wide flows and the full topology under DC power flow. Small, stylized grids and short horizons therefore misplace the location, timing, and cross-zonal propagation of the bottlenecks an ISO faces — as our coal-retirement results in Section \ref{sec: coal retirement} illustrate — which is why a realistic full-network model under projected load is necessary \cite{ARXIV251023586}. The temporal-redispatch model introduced in this paper addresses these gaps. It focuses on temporal-only flexibility as a distinct problem class motivated by single-zone transmission-level congestion; it isolates systematic flexibility from operational flexibility; and it applies the resulting program to a full annual horizon at hourly resolution to the synthetic transmission system under the ERCOT forecast. The framework is positioned not as a replacement for those more detailed physical models but as a complementary decision-support layer that uncovers the tradeoff structure between economic willingness and technical capability — a tradeoff that operators will need to quantify as data center interconnections accelerate in regions such as ERCOT.

\section{Model}
\label{sec: model}

Building on the motivation established above, we develop an operational-timescale model with a fixed grid topology and a generation fleet provisioned with sufficient capacity to meet the growing data center demand to quantify the operational economic and technical impacts of data center flexibility under explicit, hourly-resolution power-flow physics. As argued in Section \ref{sec: literature}, we restrict the framework to temporal redispatch.

Before presenting the formulation, let us clarify a point regarding the modeling perspective. The model is not that of a grid operator giving directions to specific data centers or affecting hourly data center operations. The model rather measures the overall effect of data center demand flexibility by considering the transmission grid operations at an hourly scale including the representation of hourly generation dispatch and hourly satisfaction of aggregate data center demand as estimated by the grid operator. In this way, the model assumes that the projected aggregate data center demand has some level of flexibility in shifting, potentially handled via a voluntary flexibility program. Moreover, we assume that there is a limit to this aggregate flexibility in the form of shifting cost and amount collectively achieved by data centers. This means that the implication of such flexibility is not necessarily an imposition upon a specific data center for a specific hourly operation \cite{IEEETSG2019, RSER2019109466, NI2025, IEEETSG2016}. The two economic penalties in the objective, applied to unserved demand and to the deviation from the estimated baseline (aggregate) data center demand, act as signals rather than as top-down commands: the model defers load only when the system-wide savings exceed the per-unit shifting penalty and as that penalty grows the framework decreases demand shifting \cite{CCE2023, AE2022}. The two physical limits introduced below, namely the time window over which deferred load must be recovered and the limit on how much the scheduled demand may exceed its baseline, represent the technical envelope the data centers may offer, governed by demand deferral tolerance and cooling headroom rather than capabilities the operator can extract beyond what the data centers permit. 

\begin{table}[tbp]
\centering
\caption{Notation for the mathematical model}
\label{tab: math notations}

\begin{tabular}
{p{2.2cm}p{12.7cm}}
\hline
Notation & Definition \\
\hline

\textbf{Indices} \\

$g$ & Index of a power generator \\
$k$ & Index of a branch (transmission line) \\
$n$ & Index of a bus \\
$t$ & Index of a time period (hour in the experiment) \\
$w$ & Index of a time window (daily, weekly, or biweekly in the experiment) \\

\textbf{Sets}  \\

$\mathcal{N}$ & Set of buses in the power system \\
$\mathcal{G}$ & Set of generators in the power system \\
$\mathcal{G}_n$ & Set of generators on bus $n$ \\
$\mathcal{K}$ & Set of branches in the power system \\
$\mathcal{K}_{n^+}$ & Set of branches that the starting bus is $n$ \\
$\mathcal{K}_{n^-}$ & Set of branches that the ending bus is $n$ \\
$\mathcal{T}$ & Set of time intervals \\
$\mathcal{W}$ & Set of time windows \\
$\mathcal{T}_w$ & Set of time intervals within a specific time window \\

\textbf{Parameters}  \\

$c_g$ & Operational cost for generator $g$ (\$/MWh) \\
$\mu$ & Penalty for unsatisfied demand (\$/MWh) \\
$\lambda$ & Penalty for deviation in data center demand (\$/MWh) \\
$\alpha$ & Limit of scheduled over actual data center demand of bus $n$ at period $t$ \\
$d_{nt}$ & Demand on bus $n$ at period $t$ (MW) \\ 
$d^{DC}_{nt}$ & Baseline (unshifted) data center demand on bus $n$ at period $t$ (MW) \\ 
$P_g^{max}$ & Maximum capacity of generator $g$ (MW) \\
$P_g^{min}$ & Minimum output power of generator $g$ (MW) \\
$P_k^{max}$  & Maximum thermal capacity of branch $k$ (MW) \\
$x_k$  & Reactance of branch $k$ (ohm)\\

\textbf{Variables} \\

$P_{gt}$ & The output power of generator $g$ at period $t$ (MW) \\
$P_{kt}$ & Active power flow on branch $k$ at period $t$ (MW) \\
$P'_{nt}$ & The amount of unsatisfied demand of bus $n$ at period $t$ (MW) \\
$P^{DC}_{nt}$ & Scheduled data center demand of bus $n$ at period $t$ (MW) \\
$\theta_{kt}$ & Phase angle of the branch $k$ at period $t$ (rad) \\

\hline
\end{tabular}

\end{table}

We introduce the notation used in the mathematical model in Table \ref{tab: math notations} and present the model as follows:

\begin{small}
\begin{align}
    \min \quad & \sum_{g \in \mathcal{G}} \sum_{t \in \mathcal{T}} c_g P_{gt} + \mu \sum_{n \in \mathcal{N}} \sum_{t \in \mathcal{T}} P'_{nt} + \lambda \sum_{n \in \mathcal{N}} \sum_{t \in \mathcal{T}} \ ( \ d^{DC}_{nt} - P^{DC}_{nt} \ )^+ \label{eq: obj} \\
    \mathrm{s.t.} \quad & P_g^{min} \leq P_{gt} \leq P_g^{max}  & \forall \ g \in \mathcal{G}, t \in \mathcal{T} \label{eq: gen cap} \\
    & P_{kt} = \theta_{kt} / x_k & \forall \ k \in \mathcal{K}, t \in \mathcal{T} \label{eq: dc pf} \\
    & - P_k^{max} \leq P_{kt} \leq P_k^{max} & \forall \ k \in \mathcal{K}, t \in \mathcal{T} \label{eq: trs lin cap} \\
    & \sum_{g \in \mathcal{G}_n} P_{gt} + \sum_{k \in \mathcal{K}_{n^-}} P_{kt} + P'_{nt} = d_{nt} + \sum_{k \in \mathcal{K}_{n^+}} P_{kt} + P^{DC}_{nt} & \forall \ n \in \mathcal{N}, t \in \mathcal{T} \label{eq: pf balance} \\
    & P'_{nt} \leq d_{nt} + d^{DC}_{nt} & \forall \ n \in \mathcal{N}, t \in \mathcal{T} \label{eq: uns cap} \\
    & \sum_{t \in \mathcal{T}_w} P^{DC}_{nt} = \sum_{t \in \mathcal{T}_w} d^{DC}_{nt} & \forall \ n \in \mathcal{N}, w \in \mathcal{W} \label{eq: dc balance} \\
    & P^{DC}_{nt} \leq \alpha d^{DC}_{nt} & \forall \ n \in \mathcal{N}, t \in \mathcal{T} \label{eq: dc cap}
\end{align}
\end{small}

The mathematical model is a linear program based on the direct current optimal power flow (DC-OPF) operations of the underlying grid. The objective function (\ref{eq: obj}) mathematically encapsulates the systematic flexibility of the system, comprising three distinct cost components. The first term calculates the total operational costs due to power generation. The second term applies a penalty $\mu$ to the amount of unsatisfied power demand. This penalty does not differentiate between base load and data center load, effectively simulating the value of any load that goes unsatisfied in the system. The third term introduces the cost of data center demand shifting, dictated by the penalty parameter $\lambda$. Crucially, this penalty is only calculated when the scheduled data center load is strictly lower than the data center load under its natural, unshifted schedule, represented by the positive deviation, i.e., $(x)^+=\max \{0,x\}$. $d_{nt}^{DC}$ reflects the expected operating draw from data centers at bus $n$ and time $t$ and lies below the facility's installed power and cooling capacity. This asymmetric penalty structure prevents double-counting the costs of temporal shifting. Because the model strictly enforces total energy conservation within a given time window in constraint (\ref{eq: dc balance}), any period where a data center curtails load (creating a positive deviation) mathematically necessitates another period where the data center schedules a load higher than its baseline to ``recover'' the deferred load.

Power flows are constrained by a suite of physical network constraints (\ref{eq: gen cap})-(\ref{eq: uns cap}) that ensure the power flow solutions are feasible within the DC-OPF approximation. Constraint (\ref{eq: gen cap}) sets the physical operating limits of active power output for each generator. Constraint (\ref{eq: dc pf}) is the DC power flow equation where power flowing through any transmission branch is a function of the phase angle difference across the branch and its electrical reactance. Constraint (\ref{eq: trs lin cap}) sets the maximum thermal capacity for each branch to prevent thermal overload and cascading grid failures; here the maximum thermal capacity is for the absolute power flow due to the assumption that each branch is bidirectional. Constraint (\ref{eq: pf balance}) maintains the nodal power balance in Kirchhoff’s Current Law. Constraint (\ref{eq: uns cap}) imposes the natural upper limit on the unsatisfied load amount.

While systematic flexibility defines the economic cost of maneuvering, operational flexibility defines the physical boundaries of what data centers are technically capable of executing without violating SLAs or reaching thermal limits on the data centers' cooling infrastructure. This is enforced through two distinct constraints (\ref{eq: dc balance}) and (\ref{eq: dc cap}). Constraint (\ref{eq: dc balance}) assumes that computational workloads are shift-tolerant batch processes, and are conserved within a time window. Over the set of time intervals within this window, the sum of the scheduled data center power demand must equal the sum of the original, actual data center demand. To fulfill the conservation of needed demand, the data centers eventually ``recover'' any curtailed load by running at a higher capacity in other periods within the time window. Constraint (\ref{eq: dc cap}) limits the amount of shited data center demand by defining the maximum amount of the scheduled data center load. This is done by capping the scheduled data center load relative to its actual baseline demand using the limit parameter $\alpha$. This limit governs how far data centers may run above its baseline profile during recovery hours to make up deferred load. Its admissible range is typically informed by the gap between data centers' baseline operating demand and their installed power and cooling capacity, together with the share of workload that is delay-tolerant. This gap is generally positive: data centers are routinely provisioned well above their average operating draw — for ramp-up, redundancy, and future growth — and inference and mixed cloud workloads in particular exhibit uncorrelated power peaks that leave substantial headroom for shifting and recovery of deferred loads \cite{ACM20242}. The gap is, however, workload-dependent. We therefore treat $\alpha$ as a sensitivity parameter rather than a fixed physical constant, sweeping $\alpha \in {1.5, 2.0, 2.5}$: $\alpha = 1.5$ represents a moderate, broadly attainable level of headroom for batch- and inference-oriented facilities, whereas $\alpha = 2.5$ is an optimistic limit corresponding to low-utilization sites, included to bound the value of operational flexibility rather than to assert that every facility can sustain it. The choice of values is inspired by industrial documentation \cite{NVIDIA2026} and estimates from previous works \cite{PP2025}. In all cases, the binding physical limit is the facility's thermal and cooling envelope — the rationale already encoded in constraint (\ref{eq: dc cap}) — and the availability of deferrable workload. By explicitly bounding the data center's flexibility, the framework creates a granular, parameterized tool for analyzing the effects of the amount of allowed flexibility. It strictly separates the technical capability to shift load (operational) from the cost function's financial penalties (systematic), offering a transparent mechanism to analyze temporal redispatch.

\section{Experiments}
\label{sec: experiment}

To validate the proposed temporal redispatch model, a comprehensive set of case studies is developed. This section provides a description of the experiment settings and analyzes the results, focusing on how the interplay between systematic and operational flexibility impacts grid behavior under large load growth. Finally, a special set of experiments with respect to coal generator retirement is studied, discussing the influence of localized generation retirement on the whole power grid.

\subsection{Experimental Setup}
\label{sec: experiment setting}

\begin{figure}[htbp]
\centering

\begin{subfigure}{0.50\textwidth}
    \centering
    \includegraphics[width=\linewidth]{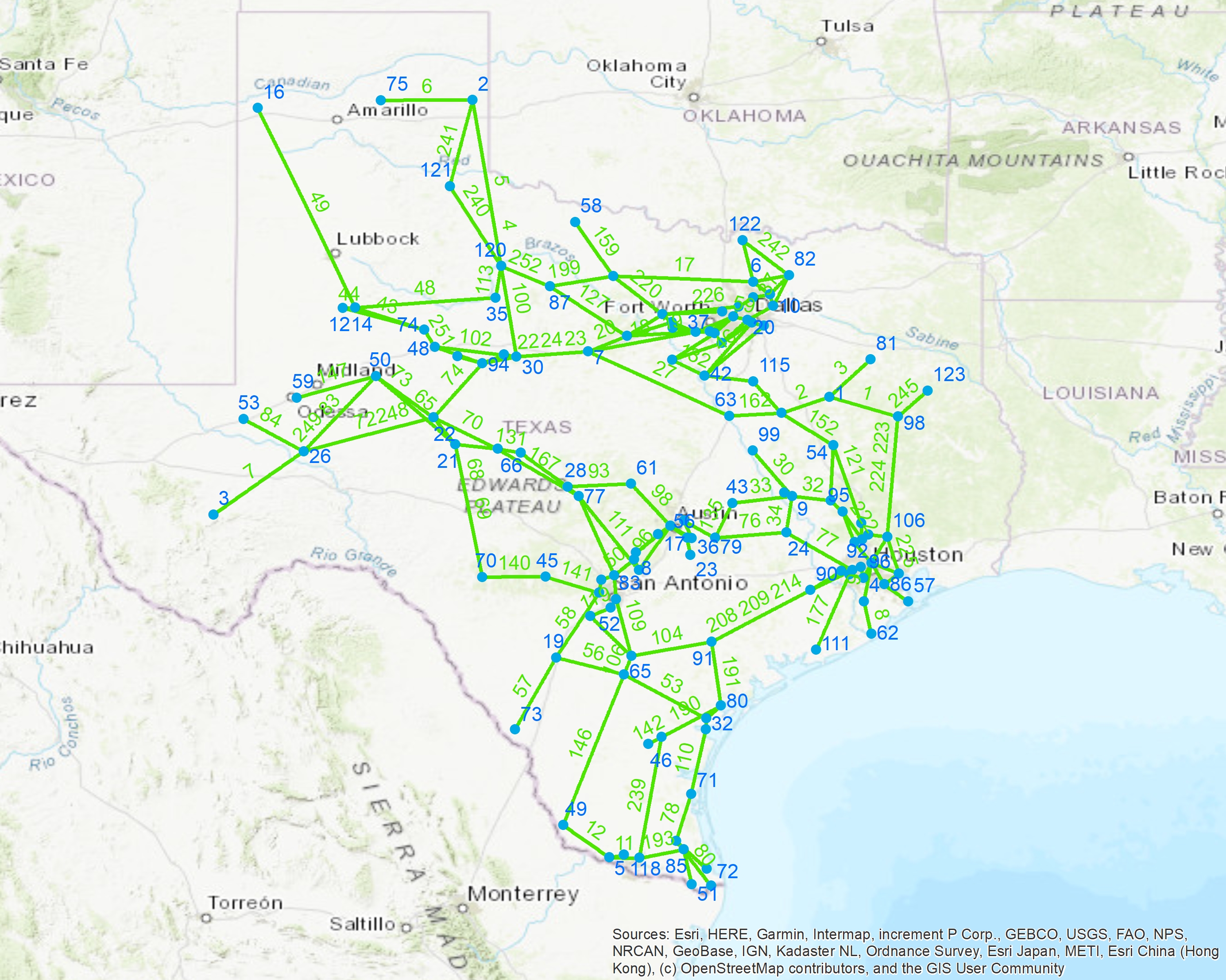}
    \caption{TX-123 BT transmission network topology \cite{SEGAN2025101774}}
    \label{fig: network topology}
\end{subfigure}
\hfill
\begin{subfigure}{0.40\textwidth}
    \centering
    \includegraphics[width=\linewidth]{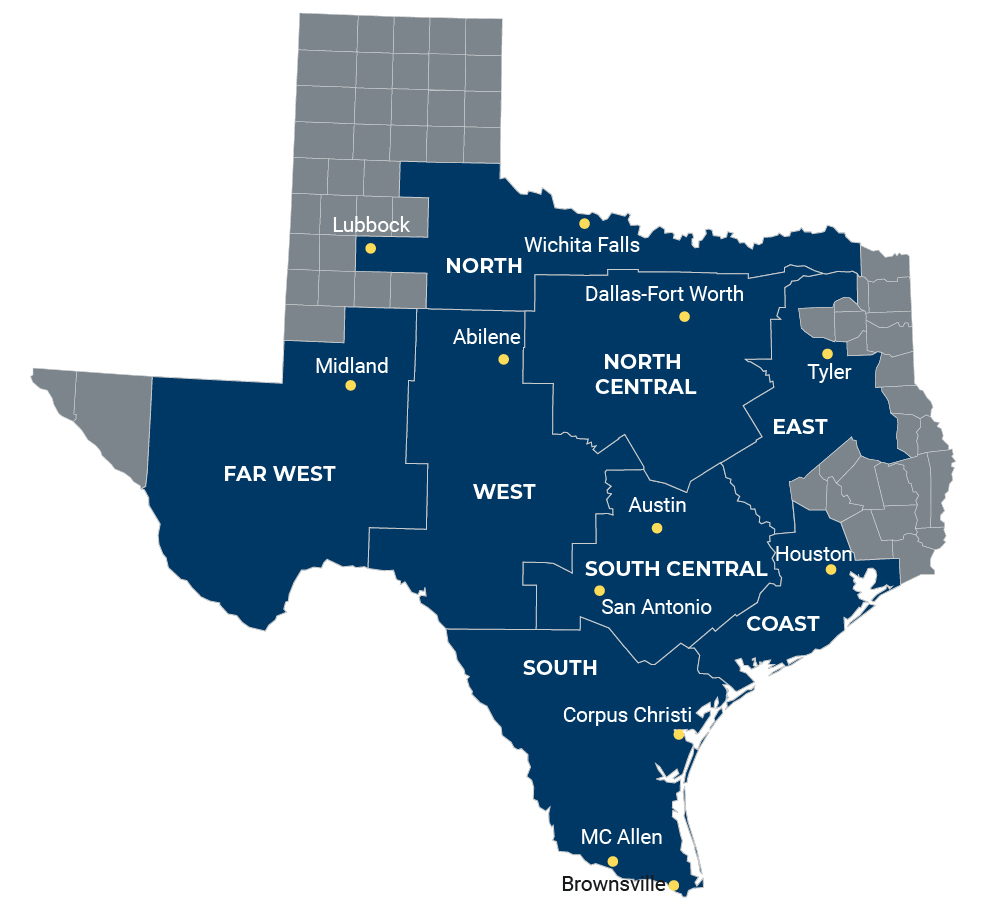}
    \caption{ERCOT weather zone topology \cite{ERCOTMAP}}
    \label{fig: zone topology}
\end{subfigure}

\caption{Experiment topology}
\label{fig: topology}

\end{figure}

The case study is constructed around a realistically constrained power grid topology, designed to mirror the localized congestion issues frequently observed in the ERCOT network. The experiments are carried out in the context of the ongoing and projected large load growth in ERCOT \cite{EIA2024}. All experiments utilize the synthetic Texas 123-bus backbone transmission (TX-123BT) system from \cite{SEGAN2025101774}, comprising 123 buses, 255 branches and 288 valid\footnote{TX-123 BT system originally contains 292 generating units, with 4 wind units coming with negligible capacity close to 0, which are omitted during preprocessing.} generating units. Figure \ref{fig: network topology} illustrates the transmission network topology, noting that multiple parallel branches may exist between a single pair of nodes.

All experiments rely on two data sources. One is the data from \cite{SEGAN2025101774}, which provides power system configuration and time series data at the \textit{transmission bus-level}, including the location and hourly load of each bus, the capacity, reactance, and terminal bus information for each branch, and the active capacity, fuel type, and operational cost of each generator. This original dataset mimics ERCOT's 2019 conditions, used as a baseline to calculate power load share of a bus within the corresponding weather zone. A key assumption here is that the load share of each bus within its corresponding weather zone remains constant. Figure \ref{fig: zone topology} provides an illustration of weather zones designated by ERCOT wherein climatological characteristics are relatively uniform \cite{ERCOTMAP}. We adopt the bus data from \cite{SEGAN2025101774} which maps each bus to a corresponding weather zone while resembling the real geographic area each load bus serves or each generator bus sits in. 

\begin{figure}[b]
    \centering
    \includegraphics[width=0.9\textwidth]{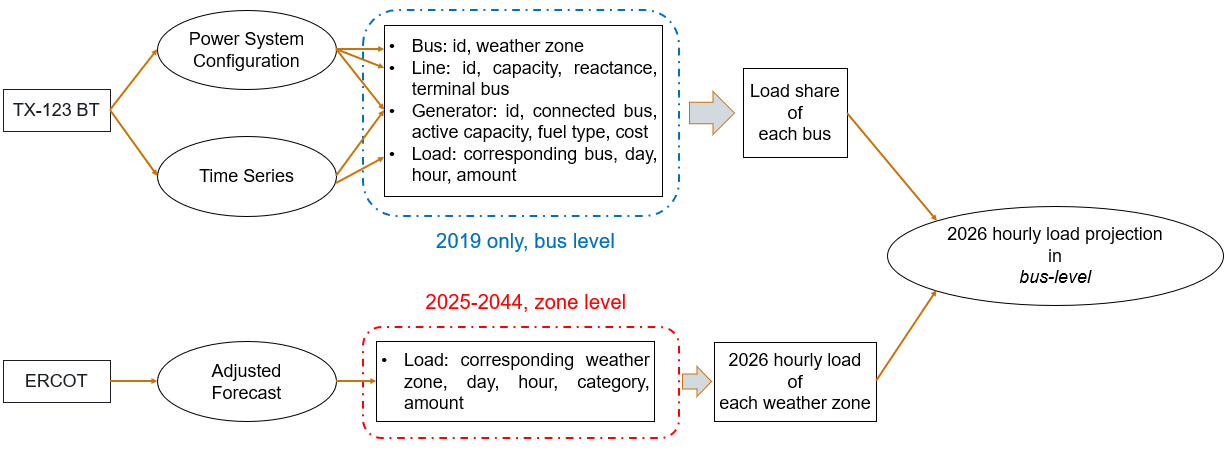}
    \caption{The workflow of bus-level hourly load forecast projection}
    \label{fig: workflow} 
\end{figure}

The second data source used for the experiments is the ERCOT Adjusted Forecast, which provides hourly load forecasts with adjustments applied to the TSP-provided large load additions with different categories at the \textit{weather zone-level} \cite{ERCOT20253, ERCOT20251, ERCOT20252}. First, we take the year 2026 as the analysis time period for all experiments. For 2026, the \textbf{Base Economic} category represents the hourly base demand of each weather zone. The \textbf{Contracts} and \textbf{Officer Letters} categories are combined to estimate the hourly additive data center demand, scaled by a factor of 0.32 based on the summer peak forecast from \cite{ERCOT20253, ERCOT20251}. This uniform scaling is appropriate because ERCOT assigns the forecast load of \textbf{Contracts} and \textbf{Officer Letters} uniformly over time within a weather zone under the same ramp schedule, negating differences in data center load share between peak and shoulder seasons. The additive part does not include the category \textbf{LFL, Large Flexible Load} because this category counts the load of existing large-load facilities, primarily cryptocurrency mining facilities \cite{EIA2024}, which is not the focus of this paper. With \textit{weather zone-level} forecast load from ERCOT and \textit{transmission bus-level} share statistics calculated from historical load provided by TX-123 BT, we are able to project the forecast hourly load at each bus. Figure \ref{fig: workflow} gives an illustration of the overall load projection workflow, and a detailed comparison between the benchmark methodology adopted in this paper to project data center load and two alternative methodologies are included in \ref{app: dc load assignment}. Following this projection methodology, we obtain load (for baseline and data centers) data at a temporal resolution of $|\mathcal{T}| = 365 \times 24 = 8760$ periods for each of the 123 buses.

To more accurately reflect 2026 grid conditions, we modify the generator and branch capacities from the power system data from TX-123BT using capacity changes by fuel type reflected in the January 2026 ERCOT data \cite{ERCOTRA}. This way, we scale the original wind and solar capacities of the TX-123BT system by factors of 1.7 and 5.4, respectively, per generator across all time periods. We assume climatological stationarity between 2019 and 2026, implying that renewable resource availability — such as solar irradiance and wind speed patterns — remains consistent. This assumption justifies the use of uniform scaling factors across all existing renewable assets rather than modeling site-specific capacity changes. By maintaining identical generation profiles, the model isolates the impacts of data center flexibility and avoids the confounding variables associated with regional infrastructure planning. To match load growth and prevent systemic flow congestion even under mild conditions, a uniform capacity expansion factor of 1.2 is applied to all existing branches. The updated TX-123BT embeds 288 generating units and a total of 162.63 GW of name-plate capacity, with 113 natural-gas units (56.00 GW), 78 wind farms (41.61 GW), 72 utility-scale solar plants (44.62 GW), 13 coal plants (14.76 GW), 10 hydro stations (0.50 GW) and two nuclear units (5.14 GW), aligning closely with 2026 ERCOT statistics.

Finally, the implementation of the demand shifting time window requires precise definition. The definition of time window builds on the time period (hour) we include in the model. For example, over an annual horizon comprising 8,760 hourly periods, a weekly time window divides the year into 52 complete weeks, starting from January 1-7, January 8-15, etc. Load shifting is strictly confined within fixed, non-overlapping time windows. Under this structure, demand at the start of a window (i.e., Hour 1, Day 1) can only be deferred to later periods within that same block; it cannot be shifted backward into a prior window. This static interval approach ensures that total energy conservation is maintained within discrete cycles, distinguishing our framework from a `sliding' or dynamic window model. For all experiments in this paper, we consider three options for the shifting time windows: daily, weekly, and biweekly. The time window for shifting is assumed to be up to 2 weeks because recovery of deferred load ideally should happen relatively soon, in a way, shifting in a day or week is preferred over a 2-week period. We believe that 2-week period is long enough for most deferred loads to find an available time slot for shifting given the daily and weekend/weekday patterns of demand and renewable generation. In a way, longer time windows will only have minimal, if any, marginal benefit of shifting flexibility at the heavy cost of making the flexibility part of the model less realistic. Hence, we think that allowing shifting of data center loads up to 2 weeks is a reasonable assumption.

\subsection{Operational Flexibility v. Systematic Flexibility}
\label{sec: of v. sf}

This section provides a comprehensive analysis of how different levels of operational flexibility impact power grid behavior under varying systematic flexibility constraints. The unsatisfied demand penalty is set sufficiently high ($\mu = 10^6$), and all instances are solved within Gurobi Optimizer version 11.0.3 on Windows 11 with a CPU model of AMD Ryzen 7 PRO 5850U. Since operational flexibility is governed by both the time window ($W$) and the scheduled data center demand limit ($\alpha$), we analyze their impacts one at a time. Specifically, the limit is fixed at $\alpha = 1.5$ during the time window evaluation, and the time window is fixed to weekly during the limit evaluation.  A supplementary sensitivity analysis on the amount of total projected data center load is provided in \ref{app: dc load amount}.

\subsubsection{Time Window vs. Systematic Flexibility}
\label{sec: time window v. sf}

\begin{figure}[htbp]
    \centering
    \includegraphics[width=0.9\textwidth]{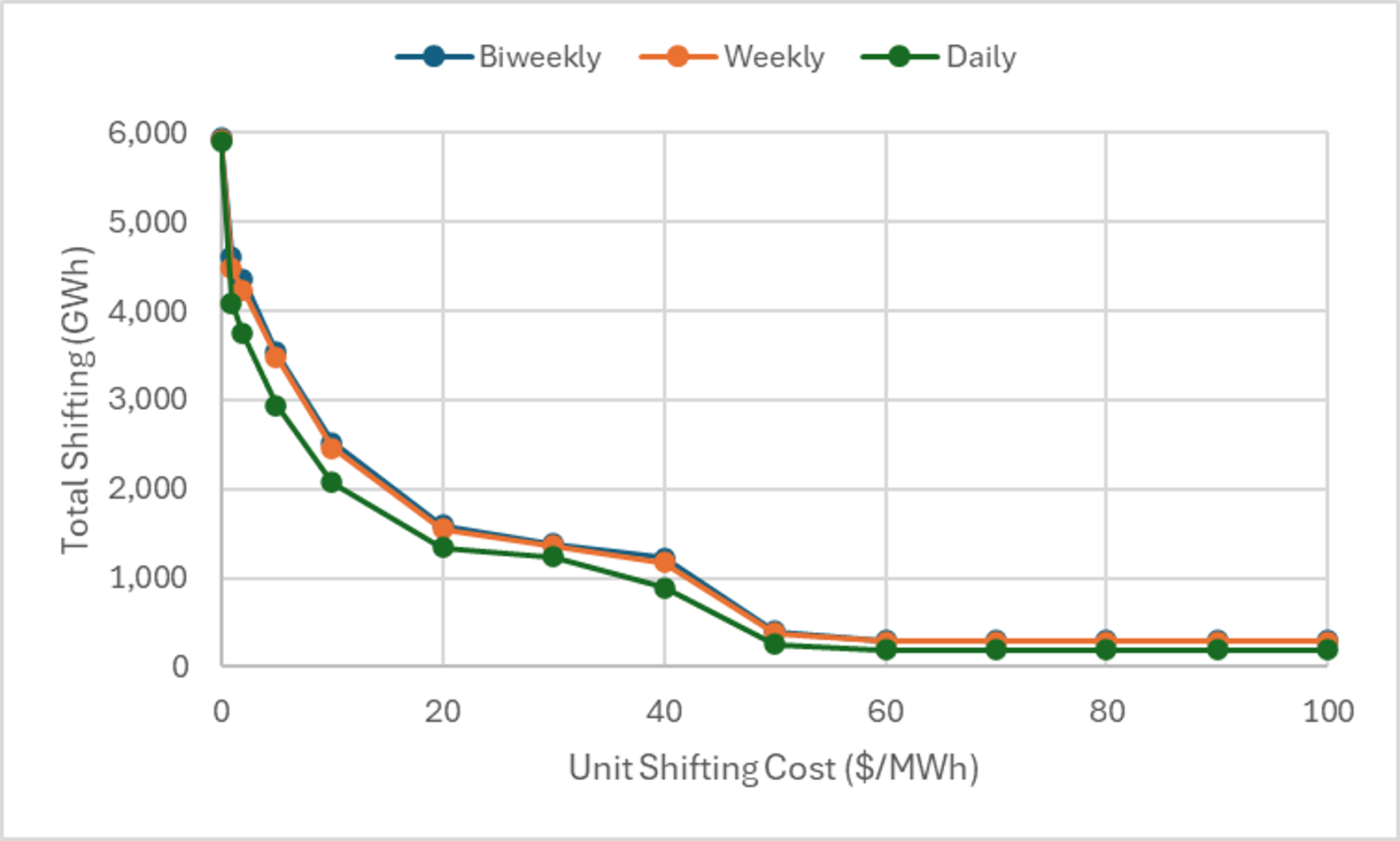}
    \caption{Total shifting amount vs. unit shifting cost over different time windows}
    \label{fig: lambda v. time window} 
\end{figure}

\begin{figure}[htbp]
    \centering
    \includegraphics[width=0.9\textwidth]{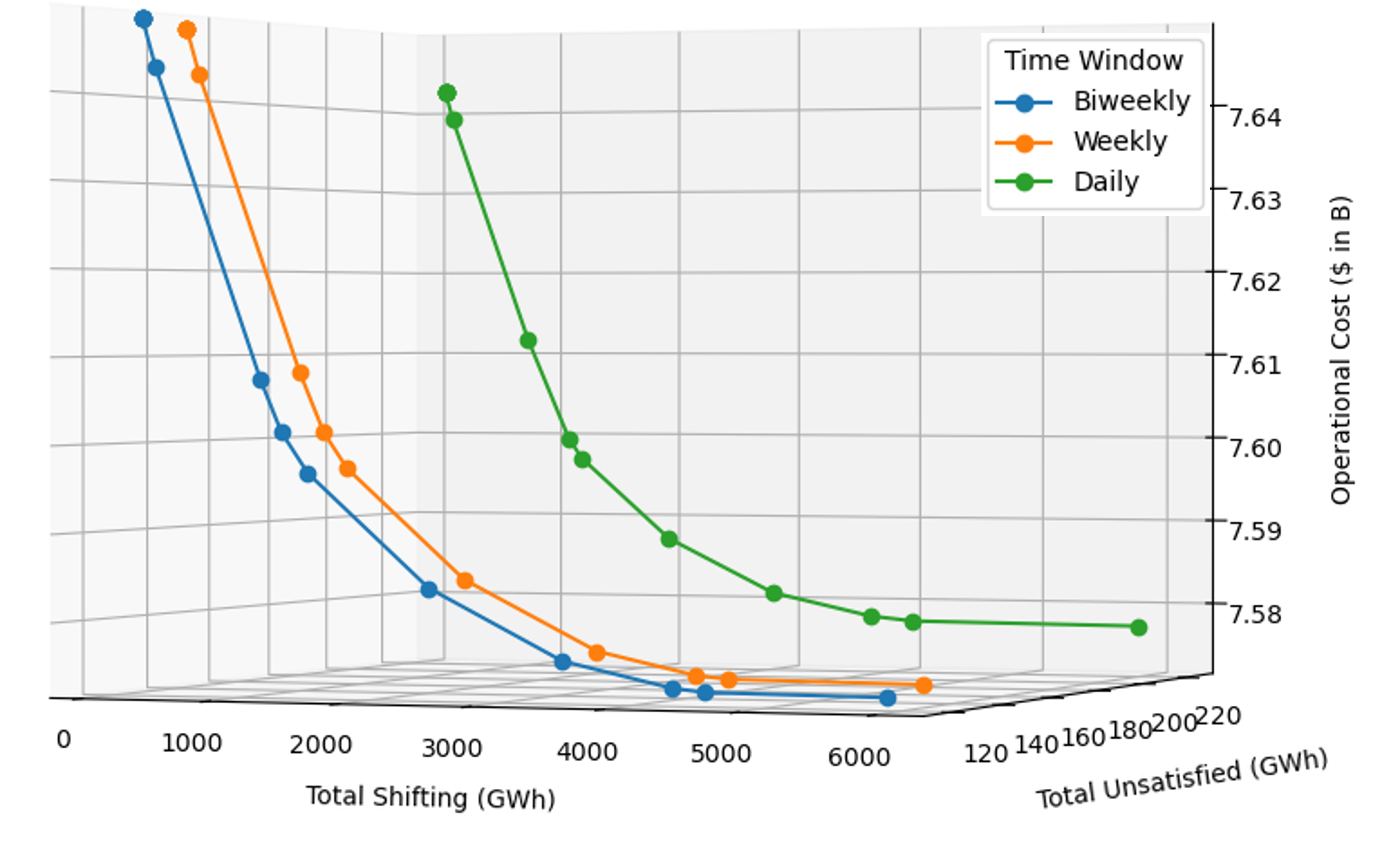}
    \caption{Operational cost vs. total shifting and unsatisfied demand over different time windows}
    \label{fig: time windows 3d} 
\end{figure}

As mentioned in Section \ref{sec: experiment setting}, three data center demand shifting time windows — daily, weekly, and biweekly — are compared under varying levels of shifting penalty ($\lambda$). Figures \ref{fig: lambda v. time window} and \ref{fig: time windows 3d} show the sensitivity of system metrics to the shifting penalty across these different time windows. Shifting penalty directly affects the amount of data center load shifting, which then affects the power generation plan and the total operational cost. Total shifting amount is calculated by summing up data center demand deviations at all buses over all time periods, i.e., $\sum_{n \in \mathcal{N}} \sum_{t \in \mathcal{T}} \ ( d^{DC}_{nt} - P^{DC}_{nt} )^+$,, and total unsatisfied demand is calculated by summing up all buses' unsatisfied demand over all time periods, i.e., $\sum_{n \in \mathcal{N}} \sum_{t \in \mathcal{T}} P'_{nt}$. Both measure the total amount over the whole year. As Figure \ref{fig: lambda v. time window} shows, the total amount of time-shifted data center demand decreases significantly as the unit shifting cost rises from 0 to small values, up to 20. The amount stabilizes (becomes nearly 0) above a certain threshold (about 50). Nevertheless, the grid utilizes load shifting with the penalty below this critical threshold, i.e., \$50/MW in our case. This is consistent for all three time windows, while a shorter time window inherently restricts the amount of shifting; the total data center load shifting in the daily window (green line) is consistently smaller than the weekly window (orange line), which in turn is smaller than the biweekly window (blue line). This is reasonable because a longer time window provides the grid with more opportunities to shift computational workloads moving power loads to more opportunistic (lower-cost generation or lower-transmission-congestion) periods.

Figure \ref{fig: time windows 3d} focuses on scenarios where the shifting penalty is below \$60/MW. As expected, longer time windows yield lower total unsatisfied demand, with the biweekly window (blue line) exhibiting the lowest levels and daily window (green line) the highest. This mirrors the trends observed in the total shifted demand amounts. Crucially, within a given window, the total unsatisfied demand remains constant regardless of the shifting penalty (and hence the total shifting). In other words, all points of the same color line are on the same straight line on the axis of total unsatisfied demand. This invariance indicates that systematic flexibility (economic incentives) cannot override the physical limitations imposed by operational flexibility constraints. No matter how inexpensive shifting becomes, the grid cannot reduce unsatisfied demand below the physical threshold dictated by the time window. Moreover, although a longer time window is able to reduce the unsatisfied demand, its influence is diminishing from daily to weekly to biweekly. Focusing on operational cost, as we move rightward along the total shifting axis for a given time window, the grid can significantly reduce total operational cost by leveraging data center demand shifting. Even then, the decreases in operational costs become smaller, as depicted by a decreasing slope of the curve as the total shift increases. 

Interestingly, as we move rightward along the axis of total shifting in Figure \ref{fig: time windows 3d}, wider time windows generally yield higher overall operational costs at the beginning (when total shifting is smaller, i.e., under 500 GWh) and then yield lower costs than shorter time windows. For example, when $\lambda = 60$ (the top point of three curves each), biweekly window has a total shifting of 296.79 GWh and an operational cost of \$7,648.30M, which is higher than the cost of weekly window \$7,647.51M (with a total shifting of 282.25 GWh) and the cost of daily window \$7,642.51M (with a total shifting of 189.38 GWh); when $\lambda = 40$ (the third top point of three curves each), biweekly window has a total shifting of 1223.29 GWh and an operational cost of \$7,607.45M, which is higher than the cost of weekly window \$7,608.36M (with a total shifting of 1173.25 GWh) and the cost of daily window \$7,611.74M (with a total shifting of 892.39 GWh). This occurs because even though a longer window offers more opportunities to shift demand, the relatively high unit penalty still restricts the solution from taking advantage of such opportunities. The optimal objective function value is achieved by reducing the unsatisfied demand, which carries the massive penalty. As we move rightward along the axis of total shifting, the unit shifting penalty decreases, and a longer time window is then able to achieve lower operational cost than shorter time windows by utilizing larger feasible regions with higher total shifting amount.

\begin{figure}[ht]
    \centering
    \includegraphics[width=0.9\textwidth]{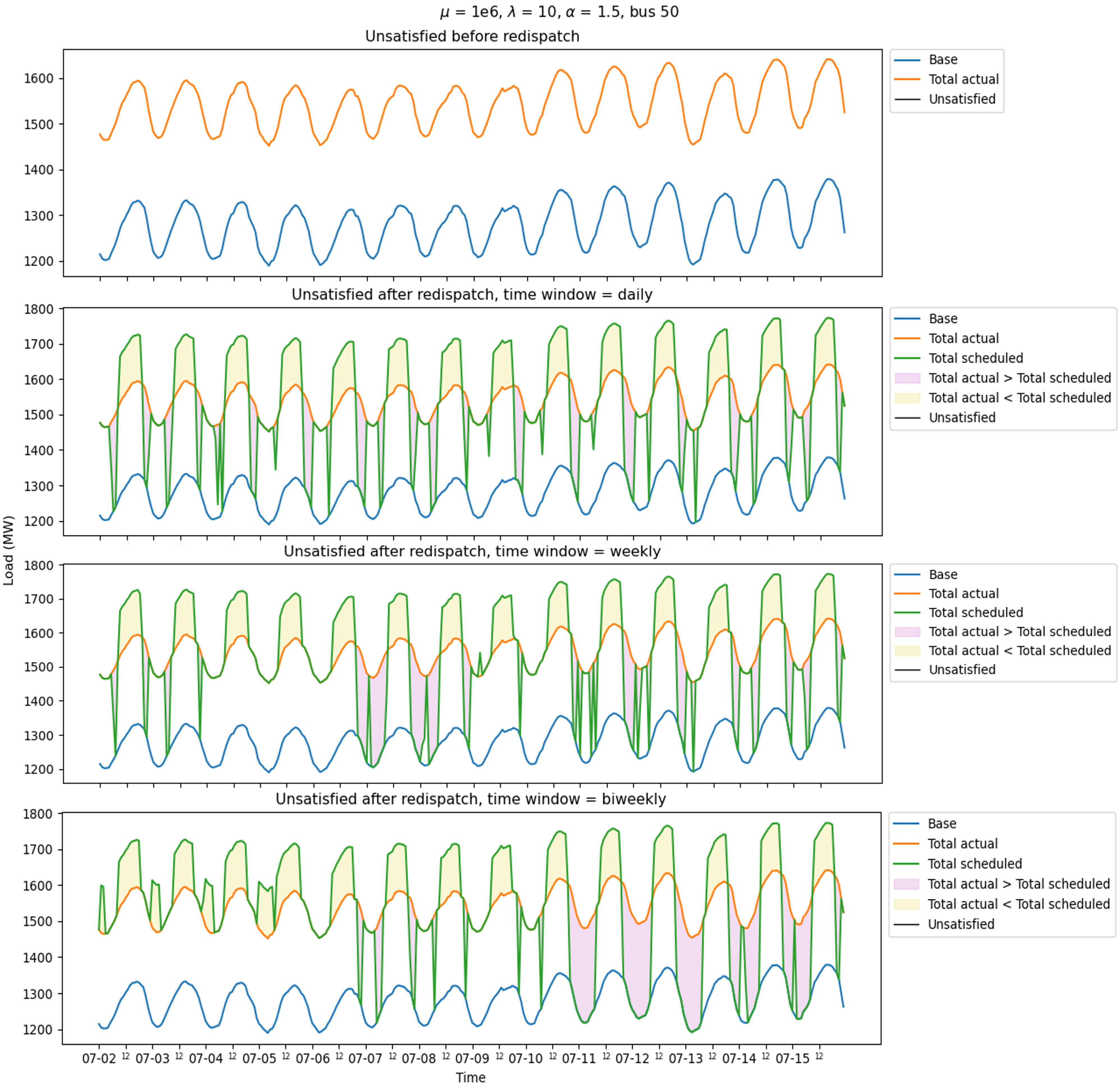}
    \caption{Base load, unsatisfied load and total load over different time windows for bus 50}
    \label{fig: time window load bus 50} 
\end{figure}

\begin{figure}[ht]
    \centering
    \includegraphics[width=0.9\textwidth]{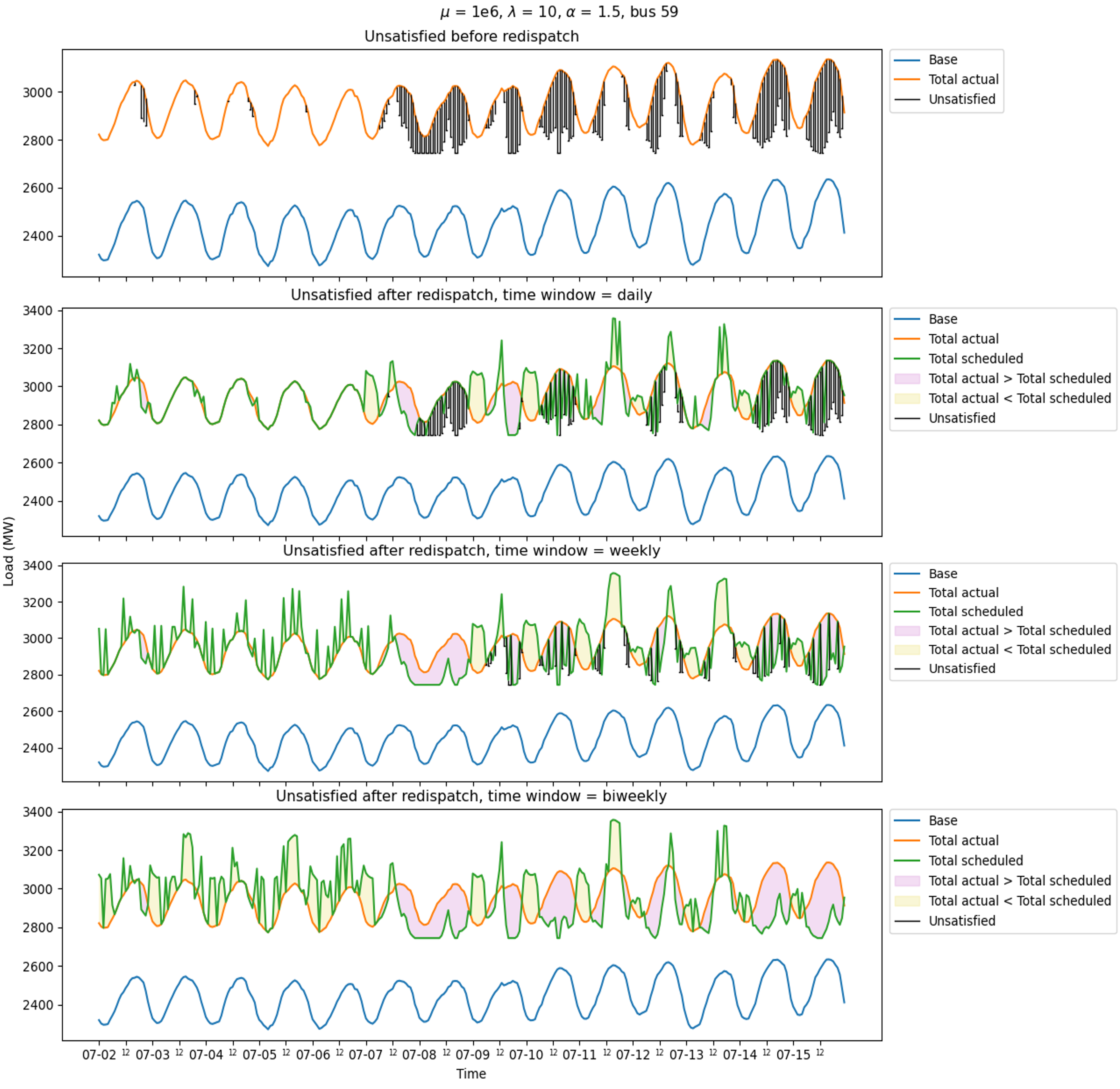}
    \caption{Base load, unsatisfied load and total load over different time windows for bus 59}
    \label{fig: time window load bus 59} 
\end{figure}

Another critical metric is the impact of data center load flexibility on the standard base demand. As Figure \ref{fig: time windows 3d} shows, the additional demand due to the introduction of data centers to the grid leads to some unsatisfied demand under the benchmark assignment. However, whether the unsatisfied demand will lead to the loss of base demand remains to be resolved. Even when there is no unsatisfied demand, how the total scheduled load curve looks like under different levels of flexibility is also an interest of analysis.  Figures \ref{fig: time window load bus 50} and \ref{fig: time window load bus 59} compare the base demand, total actual demand, total scheduled demand, and unsatisfied demand across the three time windows for two buses (Figure \ref{fig: time window load bus 50} for Bus 50 and Figure \ref{fig: time window load bus 59} for Bus 59). Bus 50 is without any unsatisfied demand and Bus 59 is with some unsatisfied demand. In fact, Bus 59 is the only bus experiencing unsatisfied demand after data center demand is introduced to the grid. The figures show bus-level data during a specific 14-day period in summer from July 2-16. This time period is a good example because it shows how the distribution of unsatisfied demand changes as we have a longer time window. The total actual demand (orange line) is the sum of base demand (blue line) and actual data center demand projected on each bus from the ERCOT forecast, and the total scheduled demand (green line) is the sum of base demand and scheduled data center demand obtained from the optimal solution. Unsatisfied demand for each hour (represented by vertical black lines) is also obtained from the optimal solution when shifting is not allowed (in the top row of each figure) or when shifting is allowed with a certain level of penalty cost ($\lambda = 10$) across the three time windows for data center demand redispatch (in the bottom three rows of each figure). We note that the optimal solution does not separate the unsatisfied demand for data center or base demand. However, we think that grid operators may use data centers' flexibility by shedding their loads first. With that in mind, when shifting is disabled, the vertical lines representing the unsatisfied demand ``start from'' the total actual demand curve. When shifting is allowed, the unsatisfied demand amounts are with respect to the total scheduled demand curve.  As Figure \ref{fig: time window load bus 59} shows, all unsatisfied demands only appear in the data center load (i.e., the black lines are within the gap between the base demand curve and the total actual or scheduled curve), meaning that it is possible to insulate standard residential and commercial demand represented in the base load from grid shortfalls. Moreover, shifting flexibility is able to reduce unsatisfied demand and a longer time window (higher operational flexibility) leads to a lower frequency and total amount of unsatisfied demand. In fact, when data center demand is flexible enough to redispatch biweekly, there is no unsatisfied demand for Bus 59 (bottom row of Figure \ref{fig: time window load bus 59}).

Beside the unsatisfied demand, the comparison between the total actual and scheduled demands also provides additional insights. Taking Bus 50 from Figure \ref{fig: time window load bus 50} as an example, when the shifting time window is set to daily, the total scheduled demand has to be matched the total actual demand for each day, i.e., each yellow area is equivalent to the sum of the pink areas for each day. Put it another way, pink areas are removed from the data centers' demand and rescheduled to yellow areas. Similarly, the sum of the yellow areas is equivalent to the sum of the pink areas for each of the weeks July 2-9 and July 10-16 when the time window is set to weekly or for the whole horizon July 2-16 when the time window is set to biweekly. As the time window is longer, the sum of all yellow (and pink) areas, i.e., the total amount of shifting, increases from 18,232.41 MWh (daily) to 18,947.93 MWh (weekly) to 19,731.34 MWh (biweekly). Although the yellow areas are similar in the peak values and appearance time across the three time windows, a longer window tends to promote pink areas to concentrate in fewer clusters with a larger area for each, usually in late evening hours. This does not mean that the solution from daily or weekly time window is infeasible in the biweekly window but the model chooses to move periodic yellow areas to larger, clustered pink areas because of the operational cost from higher operational flexibility. For a grid operator, a narrow time window such as daily will likely to result in a final (scheduled) demand profile that will have the typical (daily) repeating pattern (e.g., the green line in the daily chart - second row - in Figure \ref{fig: time window load bus 50}). A wider window will allow more ``chunky'' shifts in demand and result in total demand patterns that may not be as periodic (as daily) due to the cost minimization goal of the model with greater flexibility in the wider time window (e.g., the green line in the biweekly chart - bottom row - in Figure \ref{fig: time window load bus 50}).

\subsubsection{Limit on Scheduled Data Center Demand vs. Systematic Flexibility}
\label{sec: upper bound v. sf}

\begin{figure}[htbp]
    \centering
    \includegraphics[width=0.9\textwidth]{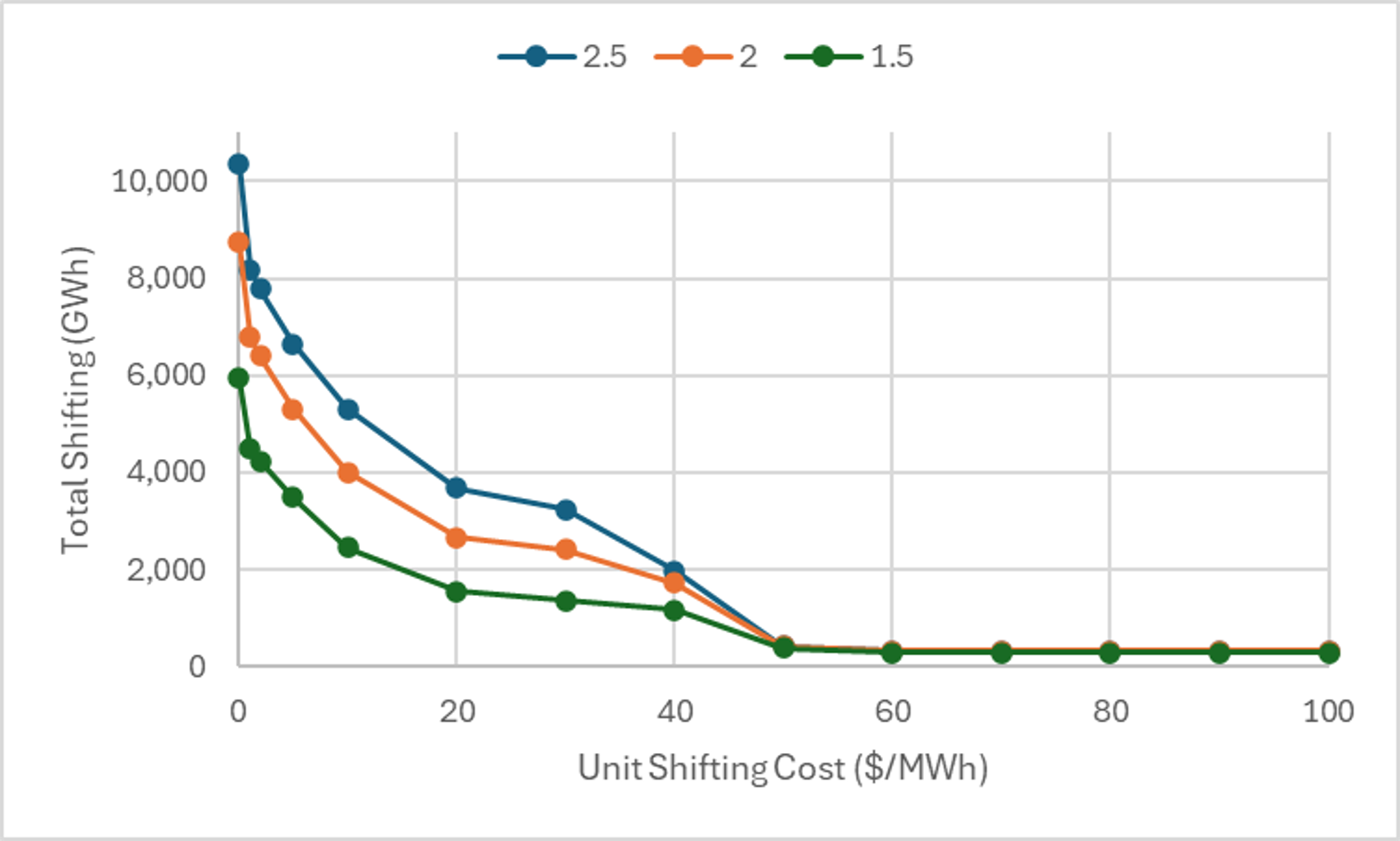}
    \caption{Total shifting amount vs. unit shifting cost over different limits ($\alpha$)}
    \label{fig: lambda v. upper bounds} 
\end{figure}

\begin{figure}[htbp]
    \centering
    \includegraphics[width=0.9\textwidth]{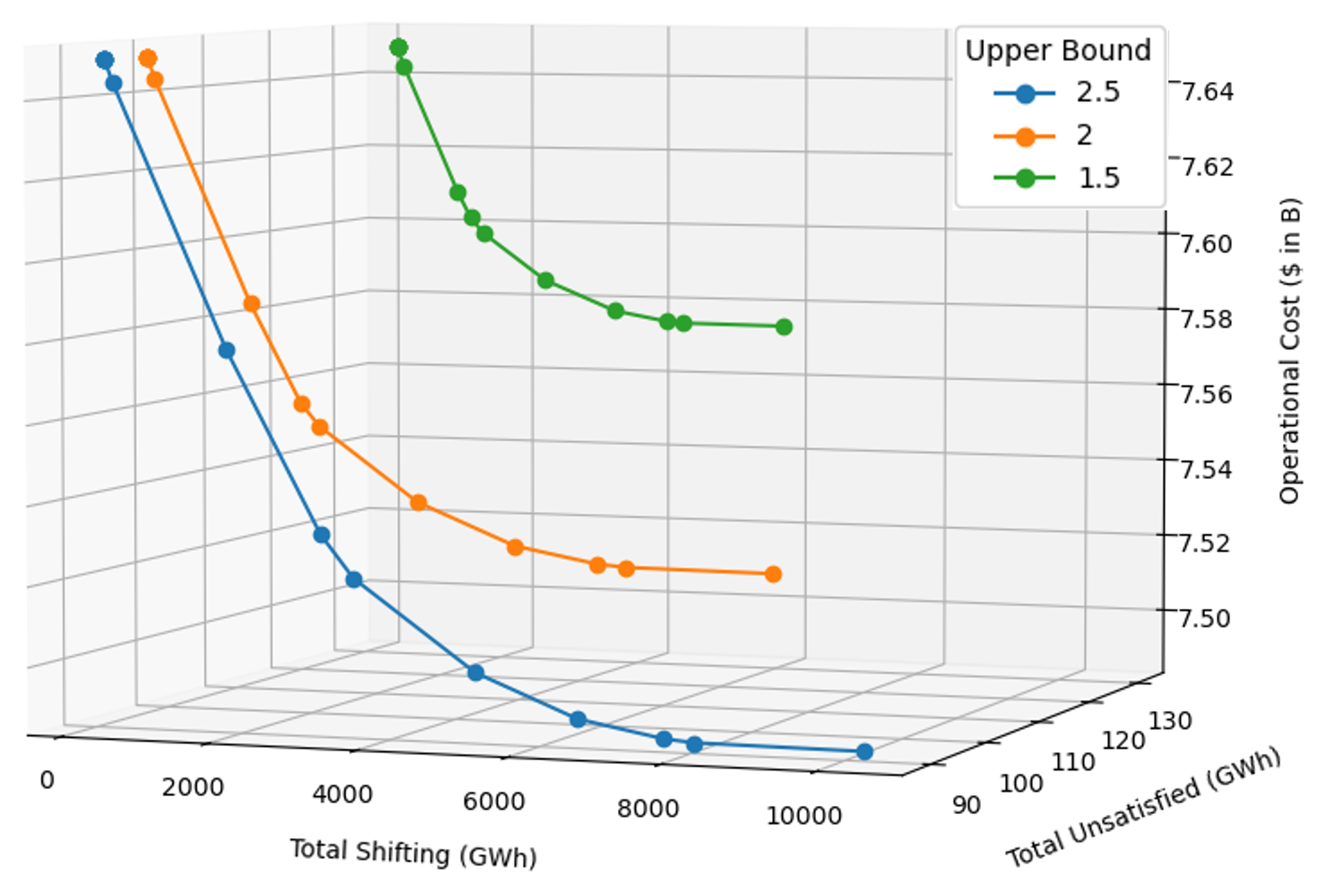}
    \caption{Operational cost vs. total shifting amount and shifting penalty over different limits ($\alpha$)}
    \label{fig: upper bounds 3d} 
\end{figure}

In this section, three levels of data center demand shifting ($\alpha$) — 1.5, 2, and 2.5 — are compared under varying shifting penalties ($\lambda$), all with the default time window of weekly redispatch. As shown in Figures \ref{fig: lambda v. upper bounds} and \ref{fig: upper bounds 3d}, the sensitivity trends are structurally similar to those observed when varying the time window. The grid utilizes load shifting when the penalty is below \$50/MW, while the total shifting amount decreases drastically as $\lambda$ increases from 0 to \$20/MWh. A conservative limit on scheduled data center demand ($\alpha = 1.5$) restricts the total shifting amount and correspondingly produces a higher total unsatisfied demand (Figure \ref{fig: upper bounds 3d}). The total unsatisfied demand for any specific $\alpha$, remains constant regardless of the shifting penalty. A higher limit on data center demand shifting ($\alpha = 2.5$) incurs slightly higher operational costs when total shifting is smaller due to higher shifting penalties ($\lambda >$ \$50/MW). However, setting $\alpha=2.5$ ultimately achieves lower total operational costs when demand shifting is inexpensive, allowing the data centers to further utilize the feasible regions. The major difference between the time window dynamics and the limit on scheduled data center demand dynamics appears in the total shifting amounts of different shifting penalties: Although shorter time windows reduce the total amount shifting, the reductions are small when compared to the differences observed as $\alpha$ becomes smaller. This makes sense because the latter controls the upper limit of shifting in each time period while the former just controls the length of time interval where shifting can occur. Within the scope of the current experiment setting, the change on time window from daily to biweekly with a restrictive scheduled data center demand ($\alpha = 1.5$) is not enough to accumulate a considerable difference in total shifting amount as it could be achieved by changing the demand limit dynamics directly.

\begin{figure}[ht]
    \centering
    \includegraphics[width=0.9\textwidth]{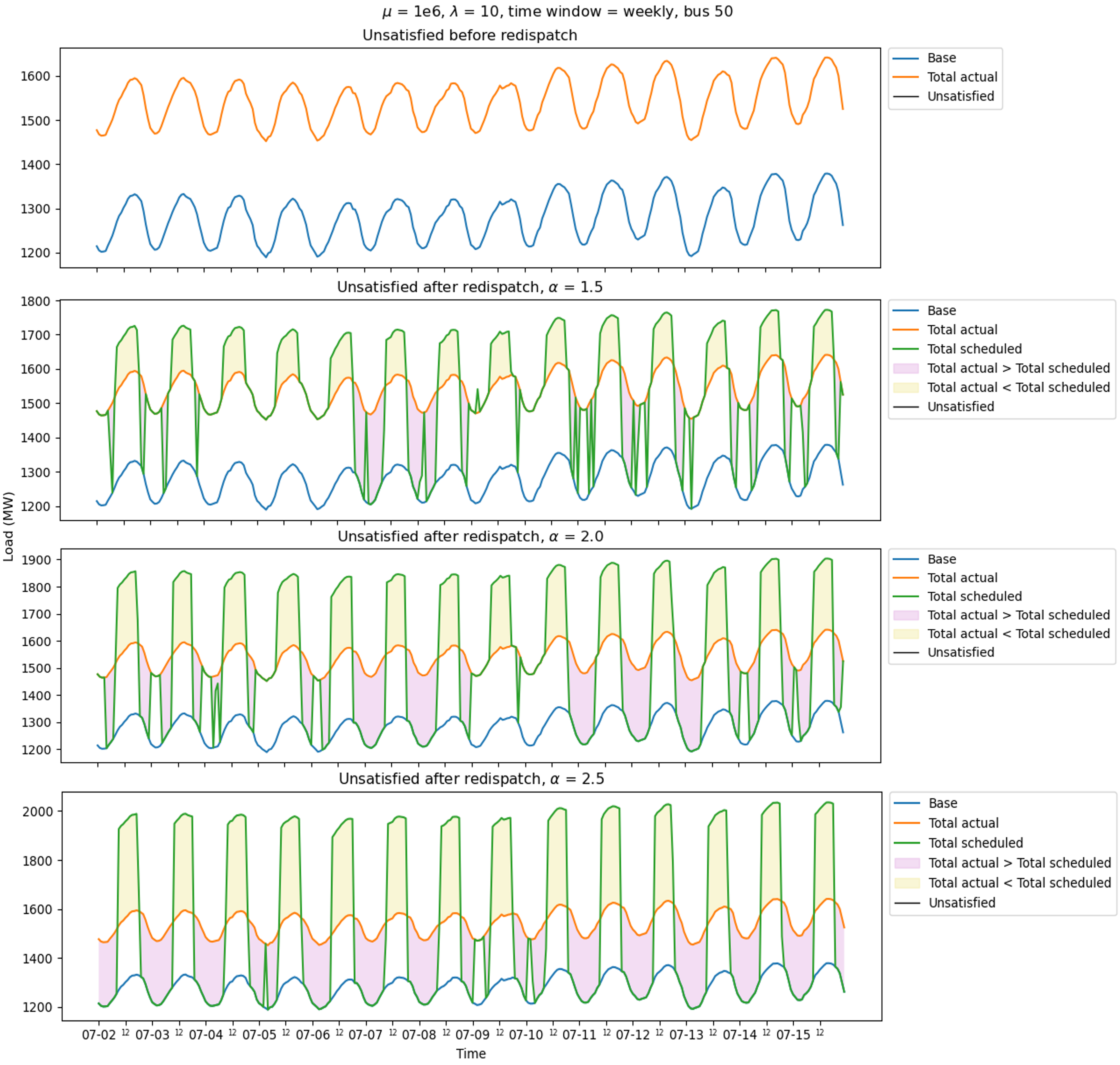}
    \caption{Base load, unsatisfied load and total load over different limits ($\alpha$) for bus 50}
    \label{fig: upper bound load bus 50} 
\end{figure}

\begin{figure}[ht]
    \centering
    \includegraphics[width=0.9\textwidth]{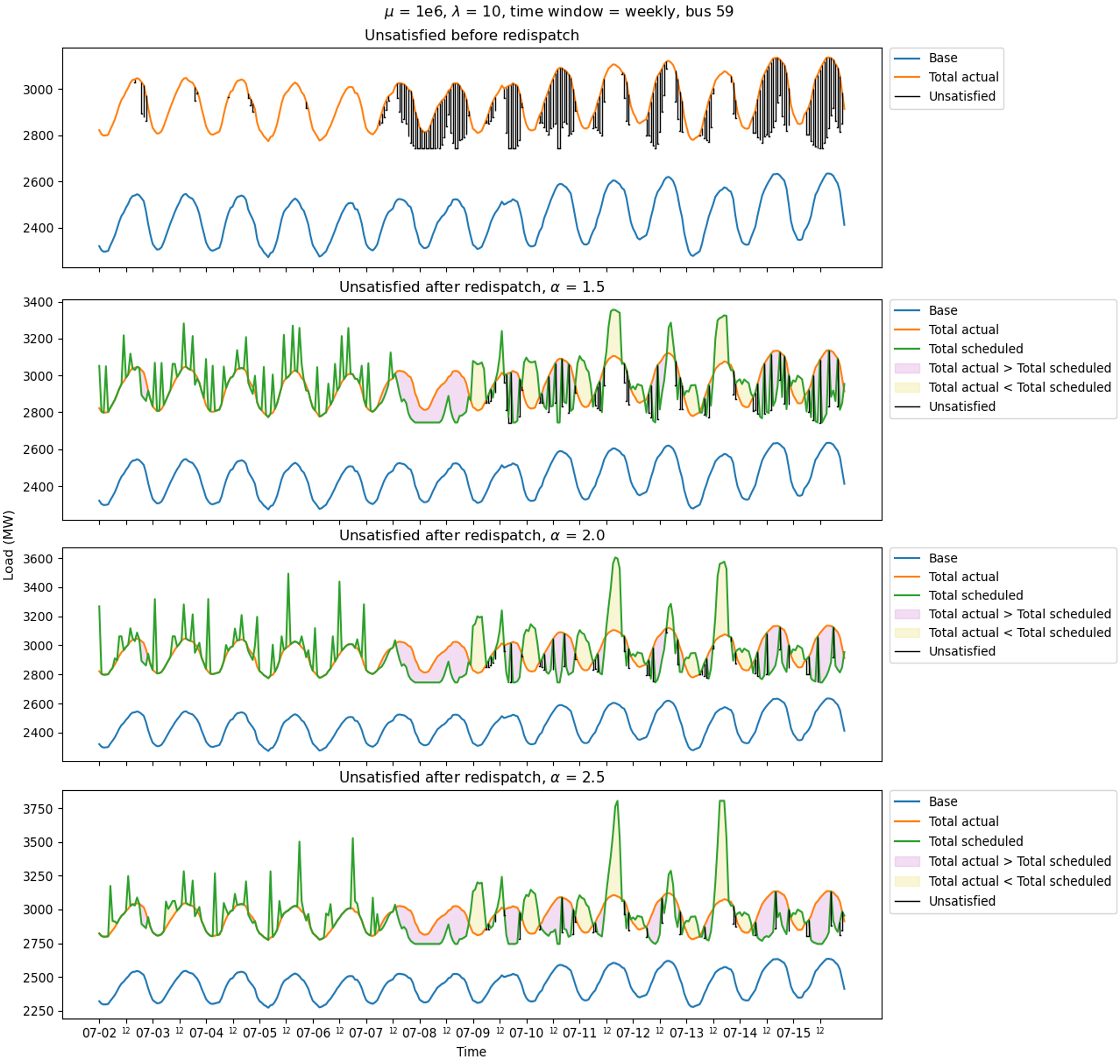}
    \caption{Base load, unsatisfied load and total load over different limits ($\alpha$) for bus 59}
    \label{fig: upper bound load bus 59} 
\end{figure}

Regarding the influence of limit of the scheduled data center demand on the power load, the results are different. Several similar patterns are still observed: first, all unsatisfied demand appears only in the data center load, meaning that it is possible to insulate the base load from grid shortfalls again; moreover, a higher operational flexibility, i.e., a larger $\alpha$, leads to a lower frequency of unsatisfied demand appearance in Figure \ref{fig: upper bound load bus 59}. The higher limit also contributes to a higher peak scheduled demand, which could be observed from the vertical axis of the bottom three rows regardless of whether an unsatisfied demand appears (Figure \ref{fig: upper bound load bus 59}) or not (Figure \ref{fig: upper bound load bus 50}). This is consistent with what we observed in Figure \ref{fig: lambda v. upper bounds} and \ref{fig: upper bounds 3d}. Nevertheless, if focusing on the total actual and total scheduled demand in Figure \ref{fig: upper bound load bus 50}, we could find that a higher operational flexibility in the limit of scheduled data center demand will promote less clustered, more periodic valley (pink) areas. This is because the grid will utilize load shifting more with a higher peak, producing larger yellow areas. Since the pink areas should be equal to yellow ares within the time window, larger yellow areas also produce more room for pink areas to grow. Although higher peaks and more valley areas may raise some concerns on whether the higher operational flexibility improves the grid stability, the periodic pattern of scheduled data center demand could still be attractive to the grid operator along with a lower operational cost.

\subsubsection{Total Net Load vs. Flexibility}
\label{sec: total net load v. flexibility}

\begin{figure}[htbp]
    \centering
    \includegraphics[width=0.9\textwidth]{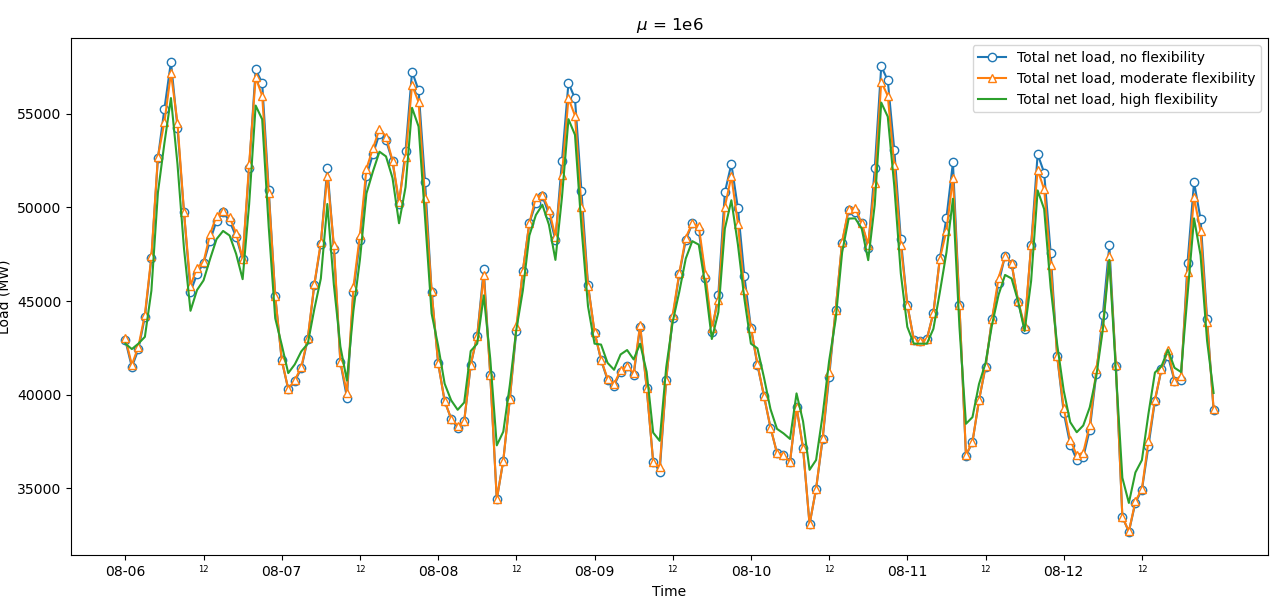}
    \caption{Total net load over different redispatch flexibilities}
    \label{fig: flexibility net load} 
\end{figure}

Both Figures \ref{fig: time window load bus 50} and \ref{fig: time window load bus 59} and Figures \ref{fig: upper bound load bus 50} and \ref{fig: upper bound load bus 59} show that higher peaks and lower valleys appear at certain buses after redispatch is enabled with either operational flexibility or systematic flexibility. On one hand, such observations are likely, since for buses in an area with rich renewable energy, such as Bus 50 and 59 in Far West, more flexibility in rescheduling data center demand will allow the system to further utilize cheap renewable energy with sufficient reserves. This is confirmed with results in Table \ref{tab: fuel utilization} in \ref{app: dc load assignment}. On the other hand, such changes to higher peaks and lower valleys can be counter-intuitive and may raise concerns about the unlikely influence of data center demand flexibility on the stability of the whole grid. That is why we provide additional details of the changes among different levels of flexibility from an aggregated, whole grid view. We consider the grid-level ``net load'' for each time period, which is defined as the total system power demand minus the power generated from renewable resources (solar and wind) for that time period \cite{EIA2014}. Figure \ref{fig: flexibility net load} compares the total net load of three flexibilities for the week Aug 6-13, one of the periods with highest peak demand in ERCOT's forecast. We choose two combinations of systematic and operational flexibilities in addition to ``no flexibility - no redispatch'' baseline for demonstration. The blue curve with circle markers depicts the net load when there is no demand shifting. The orange curve with triangle markers represents the net load when the shifting cost ($\lambda$) is \$30/MWh, the shifting time window is daily, and the limit of scheduled data center demand ($\alpha$) is 1.5, denoted as ``moderate flexibility.'' Finally, the green curve shows the net load when $\lambda$ is \$1/MWh, the time window is weekly and $\alpha$ is 2.5, called ``high flexibility.'' The results show that the reallocation of demand with higher flexibility shaves off the peaks of the total net load, within the limits of flexibility parameters, taking advantage of renewable generation when it can. When a moderate level of flexibility is enabled, the change usually occurs at lower peaks. As we increase flexibility, the differences among peaks get larger, and the differences among valleys become more visible. The figure of the total net load indicates two things: First, the higher peaks and lower valleys appear locally at buses where the area depends more on renewable resources. Second, once the variable renewable resources are excluded, higher flexibilities in redispatch promote the grid stability by flattening the net load curve with lower peaks and higher valleys. This benefit of flexibility on the power grid is in addition to the lower operational cost.

\subsection{Coal Retirement}
\label{sec: coal retirement}

The phase-out of coal-fired generation is an accelerating global trend, and the ERCOT grid is undergoing a similar transition. This section analyzes the power grid's performance under projected coal retirements when integrated with flexible data center loads. Since the generating units in the synthetic TX-123BT system do not perfectly map to actual ERCOT assets, we apply a systematic heuristic for coal retirement. For each experiment, the coal retirement only happens in one of the six weather zones with coal generating units. The zones without coal generators, Far West and West, will not be discussed in this context. If a weather zone such as North and South has only one coal generator, it will be retired, leaving no coal generators in that zone. If a weather zone has more than one coal generator, the retirement will start from the one with the largest nameplate capacity and continue to retire units until the remaining nameplate capacity is lower than the half of the original level. Table \ref{tab: retirement in each zone} provides a summary of the retirement plan for each possible weather zone used in this section's experiments. All experiments are solved within the same environment setting, under a parameter setting of high unsatisfied demand penalty ($\mu = 10^6$), free demand shifting ($\lambda = 0$), weekly conservation time window, and a limit of scheduled over actual data center demand of 1.5 ($\alpha = 1.5$).

\begin{table}[htbp]
\centering
\caption{Retirement of generators in each weather zone}
\label{tab: retirement in each zone}

\begin{tabular}{|c|c|c|}
\hline
Weather Zone & Retired Generator ID, 1st & Retired Generator ID, 2nd \\
\hline
Coast & Gen 26 & N/A \\
East & Gen 181 & N/A \\
North Central & Gen 185 & Gen 187 \\
North & Gen 269 & N/A \\
South Central & Gen 241 & Gen 261 \\
South & Gen 250 & N/A \\
\hline
\end{tabular}

\end{table}

Table \ref{tab: retirement influence in each zone} summarizes the impact of coal generator retirements on the total unsatisfied demand and influenced areas for each weather zone. Column \textbf{Weather Zone} lists the specific location where the retirement happens; column \textbf{Ori. Unsat.} and \textbf{Ori. Unsat. Area} refers to the total original unsatisfied demand over the entire grid and where it specifically happens before any retirement is applied, which is in bus 59 of the far west area; similarly, column \textbf{Cur. Unsat.} and \textbf{Cur. Unsat. Area} refers to the updated unsatisfied demand and where it happens after coal retirement; column  \textbf{Base Affected} lists whether the retirement will affect the fulfillment in base load, as the original power grid before any retirement has unsatisfied demand only for the data center load (as reported earlier). For example, for coast area, the coal retirement does not impact the unsatisfied demand in far west area, nor the base demand is influenced; nevertheless, for north area, besides the original unsatisfied data center demand from far west, the retirement also impacts the base demand for both far west and north areas, and the total unsatisfied demand increases from 132.39 GWh to 133.03 GWh. Looking into the specific topology of each weather zone, we have several interesting discoveries: for areas where the base demand remains fully satisfied after coal retirement is applied, the same node (i.e., the corresponding bus connected with a generator) may (Coast) or may not (East and South) have enough capacity from alternative fuel such as natural gas; however, for all weather zones experiencing additional unsatisfied base demand after retirements like North Central, North, and South Central, there must be a localized deficit of substitute capacity at the specific \textit{node}. This does not matter whether there is enough substitute capacity in the same \textit{weather zone}. For example, south central has enough substitute capacity from natural gas generators in the same zone from different nodes, the base load is still influenced by retirement and the influence even appears in another area, i.e., west, than the area where the retirement happens. North Central and North also show similar results that more areas are influenced by local retirement. This phenomenon underscores a critical insight for grid operators: due to transmission flow congestion, the reliability impacts of generator retirements can propagate widely, causing far-reaching disruptions that extend well beyond the immediate vicinity of the retired plant. It emphasizes that while data center flexibility can buffer the grid against general capacity shortages, it cannot entirely resolve intense, localized transmission bottlenecks exacerbated by the loss of legacy baseload generation.

\begin{table}[htbp]
\centering
\caption{Influence of coal retirement in each weather zone}
\label{tab: retirement influence in each zone}

\begin{tabular}
{|p{2.4cm}|p{1.9cm}|p{1.9cm}|p{2cm}|p{2.4cm}|p{2.4cm}|}
\hline
Weather Zone & Ori. Unsat. (GWh) & Ori. Unsat. Area & Cur. Unsat. (GWh) & Cur. Unsat. Area & Base Affected\\
\hline
Coast & 132.39 & Far West & 132.39 & Far West & No \\
\hline
East & 132.39 & Far West & 132.39 & Far West & No \\
\hline
North Central & 132.39 & Far West & 137.25 & Far West & Yes \\
 & & & & North Central & Yes \\
 & & & & North & Yes \\
 & & & & West & Yes \\
\hline
North & 132.39 & Far West & 132.45 & Far West & Yes \\
 & & & & North & Yes \\
\hline
South Central & 132.39 & Far West & 133.03 & Far West & Yes \\
 & & & & West & Yes \\
\hline
South & 132.39 & Far West & 132.39 & Far West & No \\
\hline
\end{tabular}

\end{table}

\section{Conclusion}
\label{sec: conclusion}

The unprecedented computational demands because of AI expansion place extraordinary strain on global power grids. The intense temporal mismatch between the rapid deployment of hyperscale data centers and the decade-long planning horizons required for physical grid expansion necessitates advanced operational strategies. To address this, we introduce a novel framework that strictly categorizes data center flexibility into two dimensions: systematic flexibility (dictated by economic penalties for unserved energy and schedule deviations) and operational flexibility (dictated by physical constraints, specifically the time window and the sceduled data center demand limit). By decoupling these dimensions, we provide a transparent framework for evaluating the trade-offs of temporal redispatch.

Our experimental analysis on the synthetic Texas 123-bus system shows that temporal redispatch of data center demand can improve grid performance, but its benefits depend strongly on the physical structure of flexibility. Longer conservation windows increase the amount of feasible shifting, reduce unsatisfied demand, and smooth peak-period stress by moving load toward lower-cost and less-congested hours. Increasing the scheduled data center demand limit further enlarges the feasible redispatch region and can lower operational costs when shifting is inexpensive, although it may also create higher scheduled peaks and greater aggregate generation variability. At the same time, the results show that systematic flexibility is fundamentally bounded by operational flexibility: for a fixed time window or demand cap, lowering the shifting penalty can change how much shifting the model chooses, but it cannot push unsatisfied demand below the physical threshold imposed by those constraints. Finally, the coal-retirement experiments show that although data center flexibility can help buffer general capacity shortfalls and protect base load in many cases, it cannot fully eliminate reliability problems caused by severe localized transmission bottlenecks, whose impacts may propagate far beyond the retirement zone.

While this study strictly isolates the temporal dimension of redispatch, future work will extend the model to incorporate spatial flexibility. Because digital workloads are delay-tolerant, they can be migrated across geographically dispersed data centers via global fiber optic networks. Additionally, future iterations of this framework will introduce more complex, realistic power system constraints. Transitioning from a DC-OPF to an AC-OPF model will provide critical insights into voltage stability and reactive power requirements. Furthermore, implementing detailed MILP techniques to model generator unit commitment constraints (such as minimum up/down times) and detailed energy storage architectures will yield a more comprehensive assessment of long-term grid sustainability.

\appendix
\renewcommand{\thesection}{Appendix \Alph{section}}

\section{Data Center Load Assignment}
\label{app: dc load assignment}

While the base load projection follows a straightforward distribution based on historical load share and ERCOT forecast, the projection of data center demand is subject to uncertainty since the original TX-123BT system does not natively designate data center load profiles. This section compares three distinct assignment methodologies to establish a realistic baseline.

\begin{figure}[htbp]
\centering

\begin{subfigure}{0.45\textwidth}
    \centering
    \includegraphics[width=\linewidth]{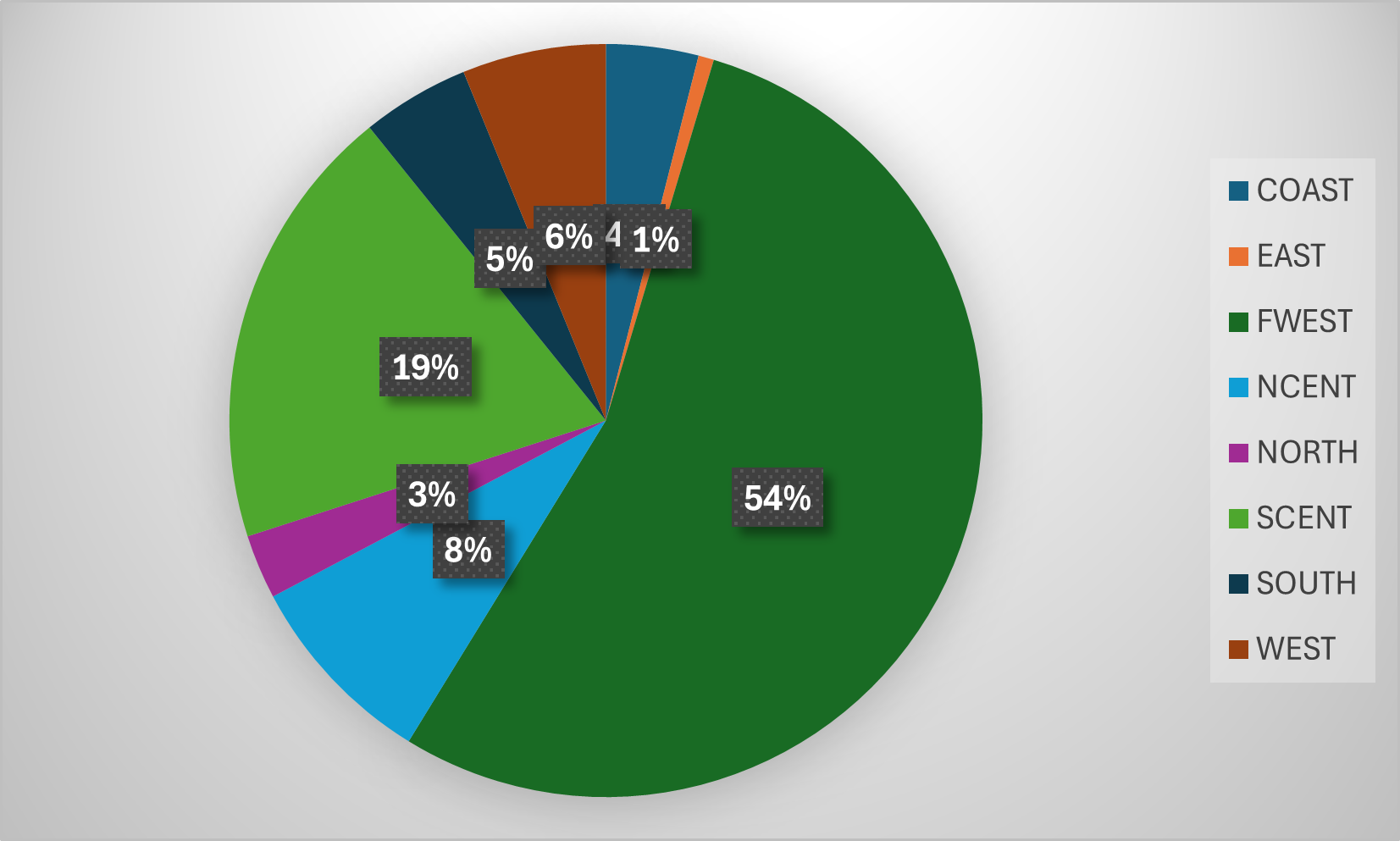}
    \caption{Data Center load share by weather zone}
    \label{fig: dc share zone}
\end{subfigure}
\hfill
\begin{subfigure}{0.45\textwidth}
    \centering
    \includegraphics[width=\linewidth]{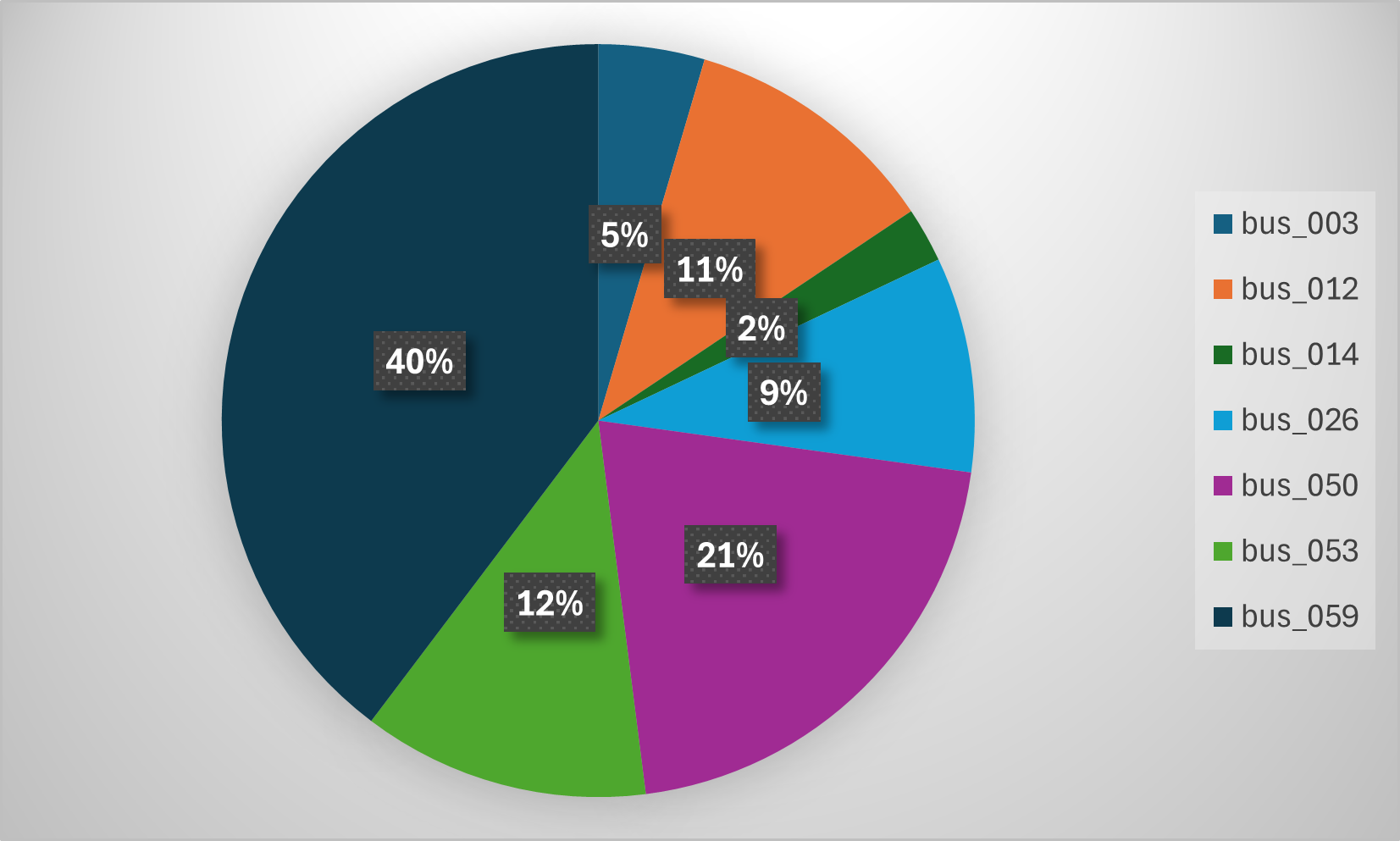}
    \caption{Base load share in Far West by bus}
    \label{fig: base share fwest}
\end{subfigure}

\caption{Load share in different spatial resolutions}
\label{fig: load share}

\end{figure}

Our benchmark methodology used in all experiments in Section \ref{sec: experiment} assigns data center demand at the bus level using the same proportional historical shares employed for base load projections. This approach ensures methodological consistency by avoiding arbitrary spatial assumptions.  Also, it can reflect the spatial concentration of large data center interconnections in the real world  \cite{ERCOT20254}, which could be supported by some statistics. As shown in Figure \ref{fig: dc share zone}, ERCOT forecasts already anticipate a spatial concentration, with the Far West zone assuming more than half of the total data center demand, higher than the second zone by more than 30\%. At higher spatial resolution, according to Figure \ref{fig: base share fwest}, historical base load shares suggest that the major demand is concentrated in only a few buses within each zone, leaving the rest buses with similar shares. Hence, allocating data center demand based on the share of base load represents a balanced benchmark.

An alternative strategy uniformly distributes data center demand across all buses within a weather zone. The assigned data center load will be purely dependent upon the total hourly load of a weather zone and the number of buses in the zone. Lacking explicit spatial data on future data center interconnections within the synthetic topology, this uniform approach avoids biased topological assumptions, though it deviates significantly from observed real-world clustering. Conversely, a third strategy allocates the data center demand exclusively to the two buses with the highest base demand shares in each weather zone, based on their normalized demand shares. Taking Far West from Figure \ref{fig: base share fwest} as an example, in this approach, only Bus 59 and Bus 50 will be assigned data center load, with a normalized share of 65.6\% and 34.4\% of the total data center load in Far West. While this approach better approximates the spatial concentration of actual data centers, it relies on heavier assumptions regarding the physical capacity of specific bus substations to handle massive point-source load injections. Here we assume that the top 2 buses with highest share have higher capacity as substations for load transmission, although the case could be really different in real world, as the synthetic topology may not be able to capture enough actual dynamics which may affect the capability of bus substations.

We evaluate these three assignment strategies assuming free data center load shifting ($\lambda = 0$) and a prohibitive penalty for unsatisfied demand ($\mu = 10^6$). The time window is set to weekly, and the scheduled data center demand cap is set to 1.5 ($\alpha = 1.5$). All instances are solved within the same environment setting. Table \ref{tab: assignment result statistics} contains the result statistics of the experiments. Branch utilization rates represent the ratio of active power flow to thermal capacity per period, while fuel utilization rates represent the average ratio of generated power to nameplate capacity across corresponding generators. High utilization period in Table \ref{tab: summary statistics} is defined as a period with branch utilization rate higher than 90\%, which means that there are at most $255 \times 24 \times 365 = 2,233,800$ periods which could be marked as high utilization periods if the corresponding power flow is over 90\% of that branch's capacity.

\begin{table}[ht]
\centering
\caption{Result statistics of each assignment}

\begin{subtable}{0.95\textwidth}
    \centering
    \begin{tabular}{|p{2cm}|p{2.5cm}|p{2.1cm}|p{2.3cm}|p{1.9cm}|p{1.8cm}|}
    \hline
    Assignment & Total Op. Cost (\$) & Total Unsat. (GWh) & Total Shifting (GWh) & High Util. Periods (\#) & Lines Involved (\#) \\
    \hline
    Equal & 7.49B & 0 & 6,373.11 & 1,275,647 & 237 \\
    Share & 7.57B & 132.39 & 5,932.76 & 1,277,625 & 237 \\
    Top 2 & 7.61B & 1300.32 & 5,488.80 & 1,278,003 & 237 \\
    \hline
    \end{tabular}
    \caption{Summary statistics}
    \label{tab: summary statistics}
\end{subtable}
\hfill
\begin{subtable}{0.95\textwidth}
    \centering
    \begin{tabular}{|p{2cm}|p{1.5cm}|p{1.5cm}|p{2.1cm}|p{1.5cm}|p{1.5cm}|p{1.5cm}|}
    \hline
    Assignment & Coal & Hydro & Natural Gas & Nuclear & Solar & Wind \\
    \hline
    Equal & 93.95\% & 99.84\% & 25.01\% & 99.41\% & 73.21\% & 96.53\% \\
    Share & 93.72\% & 99.88\% & 25.29\% & 99.37\% & 72.15\% & 96.44\% \\
    Top 2 & 93.52\% & 99.83\% & 25.36\% & 99.32\% & 71.19\% & 96.21\% \\
    \hline
    \end{tabular}
    \caption{Utilization rate by fuel}
    \label{tab: fuel utilization}
\end{subtable}

\label{tab: assignment result statistics}
\end{table}

As Table \ref{tab: fuel utilization} shows, regardless of which assignment is adopted, there are always high overall generation reserves, mainly from natural gas generators. Because the model does not take transmission cost into account, the difference on unsatisfied demand in Table \ref{tab: summary statistics} under such a huge unit penalty could only happen when the branch capacity becomes bottleneck. This is supported by the number of high utilization periods: although the number of branches involved remains the same for all three assignments — suggesting that different data center load assignments do not change the influenced lines — fewer high utilization periods are observed when data center load is allocated more uniformly. Because shifting is free while only available in temporal resolution, more buses with similar assignment means that there will be more options for model to achieve a lower total cost by avoiding penalty from unsatisfied demand and finding cheaper time periods with more buses and less congested branches. In other words, the difference on the load assignment plays a role of spatial redispatch here in advance of problem solving. These statistical results align with previous literature demonstrating that a higher spatial flexibility yields a reduction in operational costs and alleviates localized transmission congestion \cite{W34065, RSER2024114793, ARXIV251108759, AE2022}. Consequently, we adopt the proportional assignment based on base load share as our benchmark, striking an optimal balance between methodological consistency and realistic spatial distribution.

\section{Data Center Load Amount}
\label{app: dc load amount}

The projected data center load amount is one of the sources of uncertainty in terms of evaluating its flexibility in temporal shifting. Although an estimated percentage could be inferred from \cite{ERCOT20253, ERCOT20251} based on ERCOT summer peak forecast of 2026 which has already taken various factors into account, the actual data center load amount can be far from the forecast, since the scale of the future data center load is highly uncertain \cite{EPRI3002031198}. This section conducts a sensitivity analysis on the amount of projected data center load, controlled by the value of scaling factor. 

The benchmark value of the scaling factor is the one used in Section \ref{sec: experiment setting}, set at 0.32 or 32\%. As part of the sensitivity analysis, we also set the scaling factor to 41\% which is the percentage projected by ERCOT for 2027, representing the largest increase in the data center load share. Then, we vary the factor to 50\%, 60\%, ..., 100\% for the rest of sensitivity levels. We rerun the experiments of Section \ref{sec: of v. sf} at the default settings with the different percentages of the scaling factor for the data center load. Table \ref{tab: dc load share and statistics} gives a summary of each sensitivity level, the corresponding total data center load amount projected in the grid, the total operational cost and the total unsatisfied demand when the shifting is enabled or disabled. The total base load is 462,218.49 GWh, and when shifting is enabled, the time window in the table shown is weekly, the limit on scheduled demand ($\alpha$) is 1.5, and the unit shifting cost ($\lambda$) is 10.

\begin{table}[htbp]
\centering
\caption{Sensitivity levels and corresponding statistics of experiment results}
\label{tab: dc load share and statistics}

\begin{tabular}{|p{1.8cm}|p{1.8cm}|p{2.5cm}|p{2.2cm}|p{2.5cm}|p{2.3cm}|}
\hline
Sensitivity Level & DC Load Amount (GWh) & Total Op. Cost No Shifting (\$) & Total Unsat. No Shifting (GWh) & Total Op. Cost With Shifting (\$) & Total Unsat. With Shifting (GWh) \\
\hline
32\% & 19,341.41 & 7.63B & 414.58 & 7.60B & 132.39 \\
41\% & 24,781.19 & 7.80B & 922.16 & 7.77B & 492.49 \\
50\% & 30,220.96 & 7.96B & 1612.78 & 7.92B & 1073.90 \\
60\% & 36,265.15 & 8.13B & 2522.33 & 8.08B & 1945.43 \\
70\% & 42,309.34 & 8.30B & 3551.49 & 8.24B & 3038.87 \\
80\% & 48,353.53 & 8.47B & 4671.03 & 8.38B & 4243.40 \\
90\% & 54,397.72 & 8.63B & 5867.06 & 8.53B & 5539.26 \\
100\% & 60,441.92 & 8.79B & 7134.74 & 8.68B & 6838.47 \\
\hline
\end{tabular}

\end{table}

\begin{figure}[ht]
    \centering
    \includegraphics[width=0.9\textwidth]{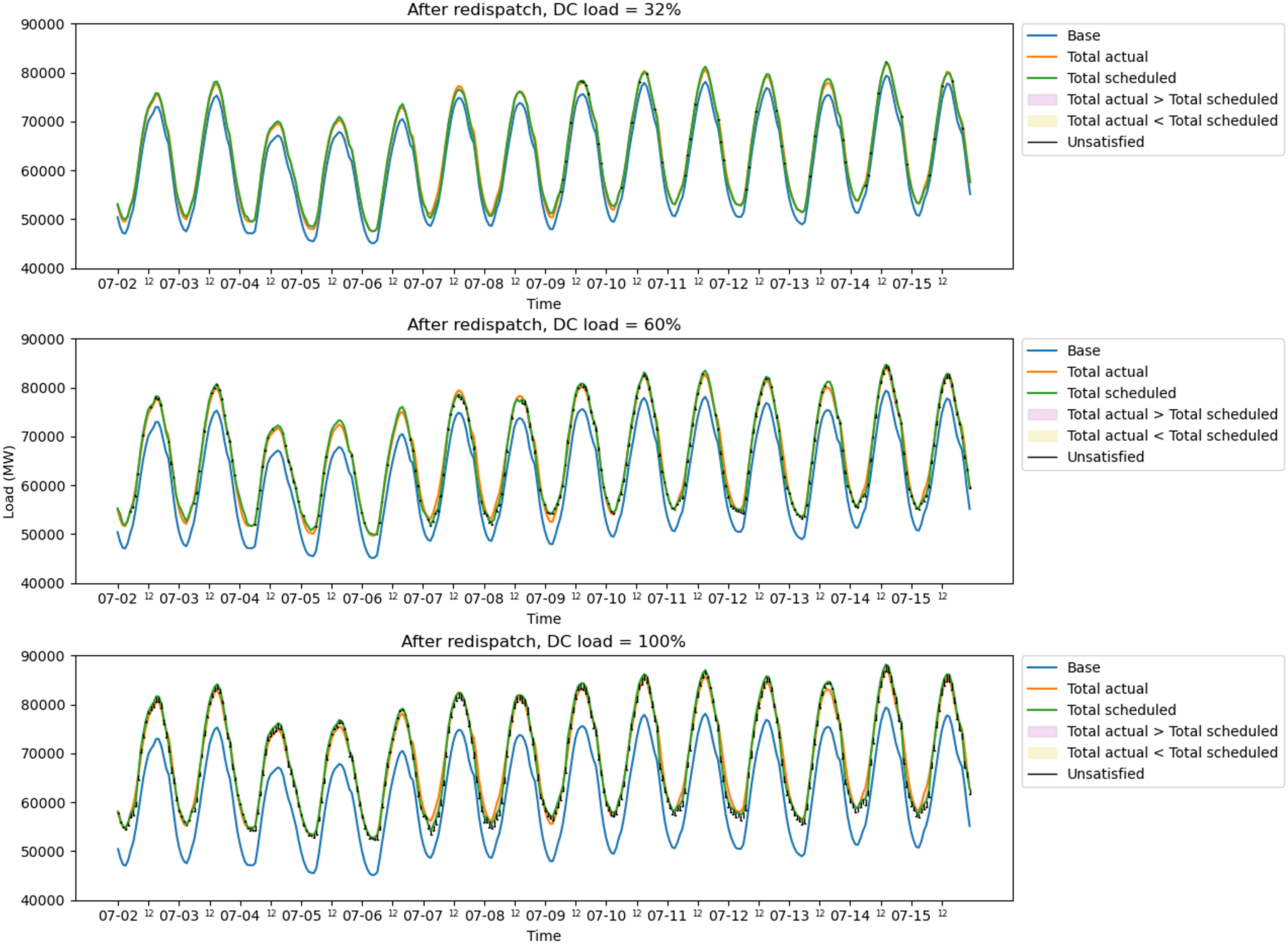}
    \caption{Total base load, total unsatisfied load and total load among different sensitivity levels}
    \label{fig: load level selected} 
\end{figure}

As the scaling factor gets larger, more data center load is projected into the system. It is expected that both total operational cost and total unsatisfied demand increase, no matter whether the shifting is enabled or not. Although more projected data center load means that the room left in the system for shifting is smaller, shifting still matters for operational cost reduction, and its influence is increasing with a larger scaling factor. If we keep increasing the projected data center load while do nothing on the generation or transmission capacity in addition to the current upgrade in the power grid, then it is reasonable to expect that shifting would play a more critical role on the cost optimization for the grid operator. Figure \ref{fig: load level selected} selects three sensitivity levels to show the total load and unsatisfied demand when shifting is enabled. Although the change in the unsatisfied demand only takes a small amount in the vertical axis of the figure, we could still clearly see its growth along with more data center load gets projected. In fact, the increase on the unsatisfied demand could threaten the satisfaction of the standard base demand on specific buses in certain cases.

\begin{figure}[ht]
    \centering
    \includegraphics[width=0.9\textwidth]{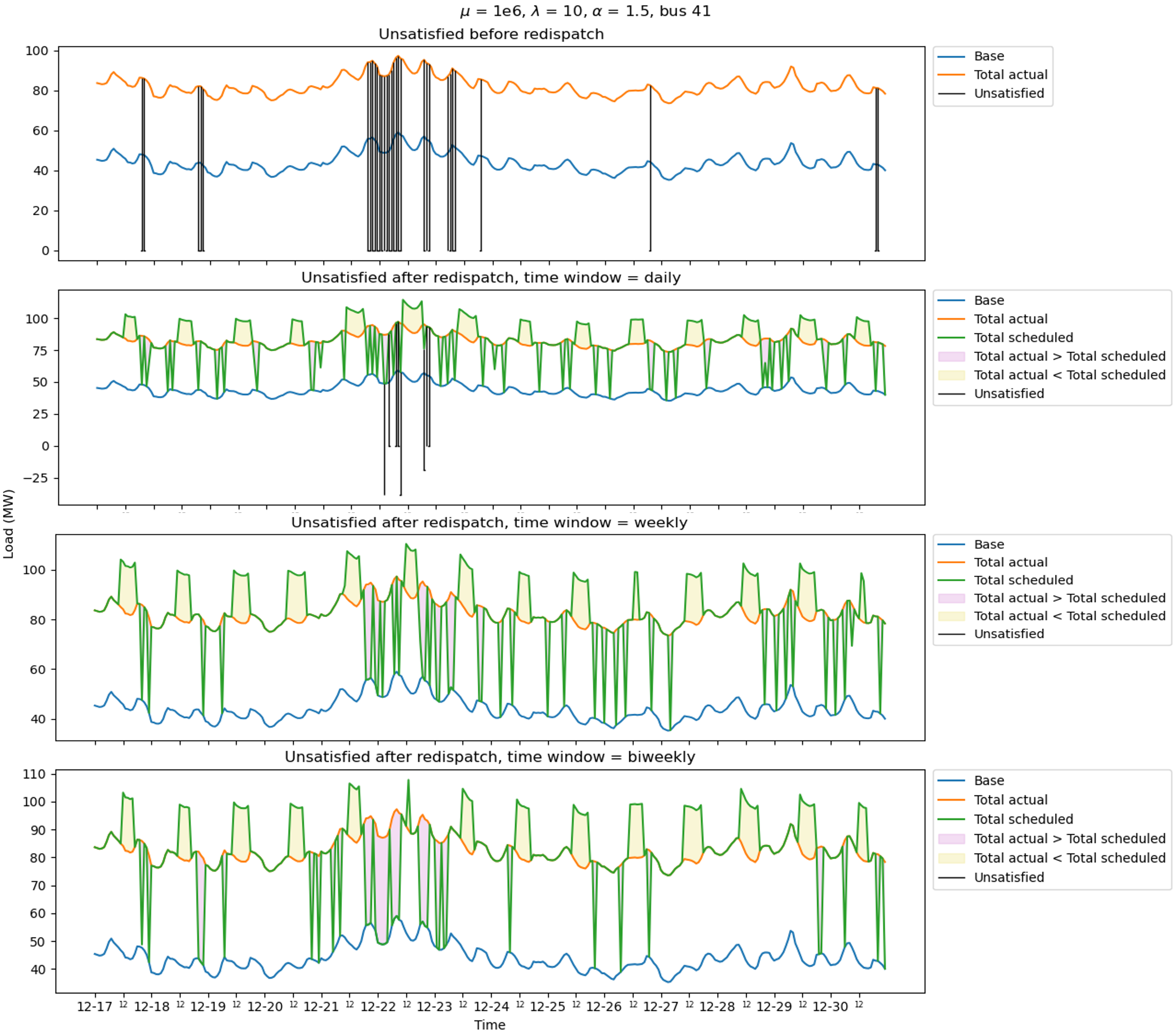}
    \caption{Base load, unsatisfied load and total load over different shifting time windows for bus 41, sensitivity level = 90\%}
    \label{fig: time window load bus 41} 
\end{figure}

Figure \ref{fig: time window load bus 41} compares the base demand, total actual demand, total scheduled demand, and unsatisfied demand across the three shifting time windows for Bus 41 in West, a bus that experiences a standard base demand plus the data center demand at 90\%. In this case, we observe that some base demand has been influenced and not been satisfied fully. Without shifting flexibility, there are days where all of the data center demand and part of the base demand is not satisfied. When the time window is set to daily, we observe a few periods with the base demand that goes unsatisfied, but the frequency of such periods is lower than when shifting is not possible. When the time window is longer, there are no unsatisfied data center or base demands. Notice that the black lines representing the amount of unsatisfied demands at several time periods go to negative. This is because the unsatisfied demand ($P'_{nt}$) at left hand side of Constraint (\ref{eq: uns cap}) is respect to the sum of base demand ($d_{nt}$) and actual data center demand ($d_{nt}^{DC}$) at right hand side, which means that $P'_{nt}$ can still be larger than $d_{nt}$ even when the scheduled data center demand ($P_{nt}^{DC}$) is 0 (at which the green line intersects with the blue line) to achieve optimality. As more data center load is projected into the grid, certain bus may schedule nothing to its own data center demand but still have a higher unsatisfied demand than its own base demand because it needs to provide power for other buses ($P_{kt}$ at right hand side of Constraint (\ref{eq: pf balance})) through branches connected with it.

With more data center load projected, another trend we observed is the further diminishing influence of operational flexibility to reduce the unsatisfied demand. In fact, when the scaling factor is at 80\% or above, the total unsatisfied demands of weekly and biweekly time windows are the same. This is consistent with what we observe from the bottom two rows of Figure \ref{fig: time window load bus 41}: the green lines follow a very similar trajectory. This is reasonable as more projected data center load leaves smaller room (fewer opportunities) for shifting. For operational flexibility in limit of scheduled data center demand, a sensitivity level of 70\% or higher will lead to the same total unsatisfied demand between $\alpha = 2.0$ and $\alpha = 2.5$, and a sensitivity level of 90\% or higher will lead to the same total unsatisfied demand among all three levels of $\alpha$. This set of analyses shows the limits of shifting flexibility: while a higher operational flexibility is often able to reduce the unsatisfied demand and that demand can be fully allocated to the data center load meaning the grid operator could insulate the base residential and commercial demand from grid shortfalls, such effect and final influence is a function of how much data center demand is expected in the grid. With a significantly large data center demand projected, the operational flexibility's benefit will be limited and the subject of the problem should turn into capacity expansion.

\clearpage
\bibliographystyle{elsarticle-num-names} 
\bibliography{references}

\end{document}